\documentclass[smallextended]{svjour3} 
\usepackage{graphicx}
\usepackage{algorithm}
\usepackage{array}
\usepackage[noend]{algpseudocode}
\usepackage{amsfonts, amsmath, amssymb}
\usepackage{bm}
\usepackage{caption}
\usepackage{citehack}
\usepackage{url}

\newcommand \D [2]{\frac{\partial #1}{\partial #2}}
\newcommand \DD[2]{\frac{\partial^2 #1}{\partial #2 ^2}}
\renewcommand{\vec}[1]{\bm{\mathrm{#1}}}
\newcommand{\V}[1]{\bm{\mathrm{#1}}}

\def \div{\nabla \cdot \mbox{}}
\def \grad{\nabla}
\def \lap{\nabla^2}
\def \half{\frac{1}{2}}
\def \ellmax{\ell_{\text{max}}}

\def \Ds{{\mathrm d}\s}
\def \Dx{{\mathrm d}\x}

\def \dt{\Delta t}

\def \dx{\Delta x}

\def \X{\vec{X}}

\def \U{\vec{U}}

\def \F{\vec{F}}

\def \I{\vec{I}}

\def \x{\vec{x}}
\def \w{\vec{w}}
\def \b{\vec{b}}
\def \s{\vec{s}}

\def \f{\vec{f}}

\def \u{\vec{u}}
\def \r{\vec{r}}
\def \e{\vec{e}}

\def \cP{\vec{\mathcal{P}}}
\def \cR{\vec{\mathcal{R}}}

\def \cL{\vec{\mathcal{L}}}

\def \cF{\vec{\mathcal{F}}}
\def \cR{\vec{\mathcal{R}}}
\def \cM{\vec{\mathcal{M}}}
\def \ctM{\widetilde{\vec{\mathcal{M}}}}
\def \chM{\widehat{\vec{\mathcal{M}}}}
\def \cG{\vec{\mathcal{G}}}
\def \cD{\vec{\mathcal{D}}}
\def \cS{\vec{\mathcal{S}}}

\def \cLIB{\vec{\mathcal{L}}_\text{IB}}
\def \ctLIB{\widetilde{\vec{\mathcal{L}}_\text{IB}}}

\def \cA{\vec{\mathcal{A}}}
\def \cAIB{\vec{\mathcal{A}}_\text{IB}}
\def \ctAIB{\widetilde{\vec{\mathcal{A}}_{IB}}}
\def \cK{\vec{\mathcal{K}}}
\def \cJ{\vec{\mathcal{J}}}

\def \xc{x_{\text{c}}}
\def \yc{y_{\text{c}}}

\begin{document}

\title{Scalable smoothing strategies for a geometric multigrid method for the immersed boundary equations \thanks{A.P.S.B.~and B.E.G.~acknowledge research funding from the National Institutes of Health (NIH award HL117063), the National Science Foundation (NSF awards ACI 1450327, CBET 1511427, and DMS 1410873), and The University of North Carolina at Chapel Hill.
M.G.K.~acknowledges research funding from the National Science Foundation (NSF award ACI 1450339) and the Department of Energy Office of Advanced Scientific Computing Research (U.S.~DOE contract DE-AC02-06CH11357).}
}

\author{Amneet Pal Singh Bhalla   \and
	    Matthew G.~Knepley          \and
	    Mark F.~Adams                  \and
	    Robert D.~Guy                   \and
	    Boyce E.~Griffith
}

\institute{Amneet Pal Singh Bhalla \at 
              Department of Mathematics, University of North Carolina, Chapel Hill, NC \\
              Carolina Center for Interdisciplinary Applied Mathematics, University of North Carolina, Chapel Hill, NC \\
              Tel.:  +1-919-962-1294\\
              Fax: +1-919-962-2568 \\
              \email{amneet@unc.edu}    
              \and
              Matthew G.~Knepley \at
              Department of Computational and Applied Mathematics, Rice University, Houston, TX \\
              \and
              Mark F.~Adams \at
              Scalable Solvers Group, Lawrence Berkeley National Laboratory, Berkeley, CA  \\
              \and
              Robert D.~Guy \at
              Department of Mathematics, University of California, Davis, CA \\
              \and 
              Boyce E.~Griffith \at
             Department of Mathematics, University of North Carolina, Chapel Hill, NC \\
              Carolina Center for Interdisciplinary Applied Mathematics, University of North Carolina, Chapel Hill, NC \\
              Department of Biomedical Engineering, University of North Carolina, Chapel Hill, NC \\
              McAllister Heart Institute, University of North Carolina, Chapel Hill, NC \\
              Tel.:  +1-919-962-1294\\
              Fax: +1-919-962-2568 \\
              \email{boyceg@unc.edu}    
}

\date{Received: date / Accepted: date}

\maketitle

\begin{abstract}
The immersed boundary (IB) method is a widely used approach to simulating fluid-structure interaction (FSI).
Although explicit versions of the IB method can suffer from severe time step size restrictions, these methods remain popular because of their simplicity and generality.
In prior work \cite{RDGuy15-gmgiib}, some of us developed a geometric multigrid preconditioner for a stable semi-implicit IB method under Stokes flow conditions; however, this solver methodology used a Vanka-type smoother 
that presented limited opportunities for parallelization.
This work extends this Stokes-IB solver methodology by developing smoothing techniques that are suitable for parallel implementation.
Specifically, we demonstrate that an additive version of the Vanka smoother can yield an effective multigrid preconditioner for the Stokes-IB equations, and we introduce an efficient Schur complement-based smoother that is also shown to be effective for the Stokes-IB equations.
We investigate the performance of these solvers for a broad range of material stiffnesses, both for Stokes flows and flows at nonzero Reynolds numbers, and for thick and thin structural models.
We show here that linear solver performance degrades with increasing Reynolds number and material stiffness, especially for thin interface cases.
Nonetheless, the proposed approaches promise to yield effective solution algorithms, especially at lower Reynolds numbers and at modest-to-high elastic stiffnesses.


\end{abstract}

\section{Introduction}
\label{sec-introduction}

Since its introduction by Peskin~\cite{Peskin72b, Peskin77} to model blood flow through heart valves, the immersed boundary (IB) method has become a widely used approach to simulate fluid-structure interaction (FSI) in a broad range of scientific and engineering applications~\cite{Peskin02}.
The flexibility of the IB approach to FSI has led to the development of many extensions, such as the ghost-cell IB method~\cite{YHTseng03}, the fictitious domain method~\cite{RGlowinski99}, the immersed finite element method~\cite{LZhang04, DBoffi08}, direct-forcing IB methods~\cite{Uhlmann05,APSBhalla13-constraint_ib}, and an IB method for immersed reactive particles~\cite{APSBhalla13-ib_reaction_diffusion}, along with other methods designed for various applications~\cite{Kim07, RMittal08, Borazjani08, APSBhalla14-EHD, DBStein16, VFlamini16-aortic_root}.
Many of these methods use 
approaches that are rooted in Peskin's original IB method, and modern versions of this method continue to see wide use, especially in biological applications.

A key feature of the IB approach to FSI is that it avoids mesh-conforming discretizations.
Instead, the IB formulation of FSI uses a single momentum equation for both the fluid and the solid, which is expressed in Eulerian form, along with a Lagrangian description of the structural deformations and resulting forces. 
The Eulerian equations are discretized on a Cartesian grid, and the Lagrangian equations are approximated on a curvilinear mesh.
Interaction between Eulerian and Lagrangian variables is mediated by discretized integral transforms with regularized delta function kernels.
These transforms \emph{interpolate} the Eulerian velocity onto the curvilinear mesh 
and \emph{spread} the structural force density to the Eulerian grid.

There has been substantial work on both explicit and implicit versions of the IB method. 
A simple version of an explicit IB time stepping scheme first uses the current configuration of the structure to evaluate the structural forces; then spreads those forces to the Cartesian grid; solves the incompressible Navier-Stokes equations; interpolates the velocities back to the structure; and finally updates the configuration of the Lagrangian mesh using the interpolated velocity field.
It is straightforward to develop more sophisticated versions of this method, e.g.~that use Runge-Kutta schemes to increase the order of accuracy of the time discretization \cite{LaiPeskin00, BEGriffith05-ib_accuracy}.
The appeal of this explicit approach is that it requires a solver only for a Cartesian grid discretization of the incompressible Navier-Stokes equations.
Fast solvers are readily available for these equations, including approaches based on fast Fourier transforms (FFT) for periodic domains and uniform Cartesian grids, or geometric multigrid (MG) algorithms for other types of boundary conditions or locally refined Cartesian grids.
Explicit time stepping is extremely effective for soft materials, but as the material stiffness $\alpha$ increases, the explicit treatment of the elastic forces imposes a time step size restriction, so that $\dt \sim \alpha^{-1}$ (for Stokes) or $\dt \sim \alpha^{-\half}$ (for Navier-Stokes). 
Although it is not straightforward to analyze the Navier-Stokes case, a simple scaling argument implies that explicit IB method in Stokes flow conditions requires $\dt \sim \dx$ for thin elastic membranes and $\dt \sim \dx^3$ for thin beams.
The stability restrictions for thick structures are less severe.

The alternatives to explicit IB methods are fully implicit and semi-implicit IB methods.
Fully implicit IB methods can allow for the stable use of any time step size \cite{Newren07}. 
It is also possible to develop stable semi-implicit IB methods that use spreading and interpolation operators defined with respect to the current structural configuration, or an estimate of the new position \cite{Newren07}.
These are referred to as \textit{lagged} IB coupling operators~\cite{MayoPeskin93, HDCeniceros09}.
In effect, this approach linearizes the geometrical nonlinearities associated with the coupling operators.
The structural configuration used to evaluate the Lagrangian forces still must be treated implicitly in such discretizations to maintain energy stability \cite{Newren07}, but the resulting system of equations is substantially simpler than that of a fully implicit formulation.

Work on implicit IB formulations dates back to the first IB methods \cite{Peskin72b, Peskin77}, 
but here we briefly review research over only the past decade to develop efficient \mbox{(semi-)implicit} IB methods.
Hou and Shi~\cite{TYHou08a, TYHou08b} proposed a semi-implicit and unconditionally stable discretization of the IB equations for steady and unsteady Stokes flow for simple periodic interfaces with linear elasticity.
They deploy a small-scale decomposition to obtain a formulation that can be expressed explicitly using Fourier transforms, which allows them to obtain an efficient solution method.
Linear solvers based on semi-implicit discretizations of the IB method have also been proposed to treat more general structural geometries.
Two notable examples are the works of Mori and Peskin~\cite{MoriPeskin08} and Ceniceros et al.~\cite{HDCeniceros09}.
Both of these studies reformulate the IB equations by eliminating the Eulerian variables, so that the systems to be solved involve only Lagrangian degrees of freedom.
An update to the Eulerian velocity and pressure is made thereafter by using the new position of the immersed structure.
Mori and Peskin~\cite{MoriPeskin08} suggest a simple diagonal preconditioner for the unstructured Lagrangian system, whereas Ceniceros et al.~\cite{HDCeniceros09} employ an algebraic multigrid solver on the unstructured Lagrangian mesh, in which coarser and finer Lagrangian meshes are obtained by adding and removing Lagrangian points from a base mesh.
Ceniceros et al.~also advocate  precomputing an explicit matrix-based representation of the Lagrangian linear operator for modest ratios of Lagrangian to Eulerian degrees of freedom.
For periodic domains, they are able to apply the Lagrangian operator efficiently by exploiting the approximate translational invariance of Peskin's regularized delta functions~\cite{Peskin02}.
The idea of precomputing the Lagrangian matrix operator has also recently been employed by Kallemov et al.~\cite{BKallemov16-RigidIBAMR} and Usabiaga et al.~\cite{FBalboaUsabiagaXX-rigid_multiblob} for rigid-body IB methods.

Although solving implicit or semi-implicit IB formulations using only Lagrangian variables can be very efficient for certain problems, developing scalable general-purpose algorithms for these formulations is difficult.
In particular, constructing multigrid methods for Lagrangian formulations of the IB equations is challenging because the systems to be solved fundamentally involve the solution operator for the Stokes equations, which is used in this formulation to eliminate the Eulerian velocity and pressure variables.
To avoid this difficulty, Guy et al.~\cite{RDGuy15-gmgiib, RDGuy12} and Zhang et al.~\cite{QZhang14} proposed multigrid preconditioners for semi-implicit IB formulations in which the Lagrangian variables are eliminated.
This approach requires the solution of Stokes-like systems of equations on structured Cartesian grids that involve only the Eulerian variables, thereby facilitating the development of geometric multigrid algorithms.
Specifically, Guy et al.~\cite{RDGuy15-gmgiib, RDGuy12} developed a geometric multigrid method for this Eulerian IB formulation similar to Vanka's method for the Stokes equations~\cite{SPVanka86}, whereas Zhang et al.~\cite{QZhang14} proposed an approximate block-factorization preconditioner for this system that is similar to block multigrid preconditioners for the Navier-Stokes equations~\cite{ElmanEtAl08, BEGriffith09-efficient, MCai14-variable_coefficient_stokes}.
The approach of Zhang et al.~does not appear to provide a robust semi-implicit solution strategy, and multiplicative smoothing strategies like that developed by Guy et al.~present limited opportunities for parallelization.
Further, Vanka-like smoothing for the IB equations requires the solution of relatively large block systems, resulting in a computational cost much greater than multigrid algorithms that can rely on simpler point relaxation smoothers, such as Jacobi or Gauss-Seidel smoothing.

This work extends the multigrid approach of Guy et al.~\cite{RDGuy15-gmgiib} 
by introducing two different smoothing approaches that are amenable to large-scale parallelization.
One smoother that we consider is similar to the Vanka-like scheme developed by Guy et al., but it uses a restricted additive Schwarz (RAS) method~\cite{XCCai99, EEfstathiou03-RAS} to couple the ``big box'' solves \cite{RDGuy15-gmgiib} required by this smoothing algorithm instead of the 
multiplicative algorithm developed by Guy et al.
This approach allows each of the subdomain solves to be processed independently.
We also develop a Schur complement (SC) smoother for the Stokes-IB equations based on an approximate block factorization.
This approach is similar to the method of Zhang et al.~\cite{QZhang14}, except that here we use the block factorization as a smoother.
We show that this SC smoother can be effective even when using only lightweight subdomain solves involving a few iterations of Chebyshev-accelerated Gauss-Seidel applied to Poisson-like operators.

Unlike the work of Guy et al.~\cite{RDGuy15-gmgiib}, here we consider nonzero Reynolds number flows in addition to the Stokes flow regime.
As in earlier studies, solver convergence rates are shown to degrade with increasing material stiffness.
This study also reveals, for the first time, that the linear solver convergence rates degrade with increasing Reynolds numbers.
For low Reynolds numbers or Stokes flows, however, mesh-refinement studies demonstrate essentially optimal scaling under only a mild CFL-type time step size restriction for a range of material stiffnesses.

\section{Immersed boundary method}

\subsection{Continuous equations of motion}

In the immersed boundary (IB) formulation of fluid-structure interaction (FSI) problems, an Eulerian description is used for the momentum equation and divergence-free condition of both the fluid and the structure, and a Lagrangian description is used for the structural deformations and the resulting structural forces.
We denote by $\x = (x_1, \ldots, x_d) \in \Omega$ fixed Cartesian coordinates, in which $\Omega \subset \mathbb{R}^d$ is the fixed domain occupied by the entire fluid-structure system in $d$ spatial dimensions.
We denote by $\s = (s_1, \ldots s_d) \in U$ the fixed material coordinate system attached to the structure, in which $U \subset \mathbb{R}^d$ is the Lagrangian curvilinear coordinate domain.
The position of the immersed structure at time $t$ is denoted $\X (\s,t) \in \Omega$.
To simplify the implementation, we consider only thin (codimension-1) massless structures and thick (codimension-0) neutrally buoyant bodies.
In the case of a thick immersed body, the fluid and structure share the same uniform mass density $\rho$, and we further assume that the structure is viscoelastic with the same dynamic viscosity $\mu$ as the fluid.
The equations of motion of the coupled fluid-structure system are~\cite{Peskin02}
\begin{align}
\rho\left(\D{\u}{t}(\x,t) + \u(\x,t) \cdot \grad \u(\x,t) \right) &= -\grad p(\x,t) + \mu \lap \u(\x,t) + \f(\x,t), \label{eqn-momentum} \\
 \div \u(\x,t) &= \V 0, \label{eqn-continuity} \\
 \f(\x,t) &= \int_{\Omega} \F(\s,t) \, \delta(\x - \X(\s,t)) \, \Ds, \label{eqn-F-f} \\
 \D{\X}{t} (\s,t) &= \int_{U} \u(\x,t) \, \delta(\x - \X(\s,t)) \, \Dx,\label{eqn-u-interpolation} \\
 \F(\s,t) &= \cF[\X(\cdot,t)](\s,t). \label{eqn-elastForce}
\end{align}
Eqs.~\eqref{eqn-momentum} and \eqref{eqn-continuity} are the incompressible Navier-Stokes equations written in Eulerian form, in which $\u(\x,t)$ is the velocity, $p(\x,t)$ is the pressure, and $\f(\x,t)$ is the elastic force density.
Eq.~\eqref{eqn-elastForce} determines the Lagrangian structural force density from the configuration of the immersed structure via a functional $\cF: \X \mapsto \F$.
Interactions between Lagrangian and Eulerian quantities in Eqs.~\eqref{eqn-F-f} and~\eqref{eqn-u-interpolation} are mediated by integral equations with Dirac delta function kernels, in which the $d$-dimensional delta function is $\delta(\x) = \Pi_{i=1}^{d}\delta(x_i)$.
Eq.~\eqref{eqn-F-f} converts the Lagrangian force density $\F(\s,t)$ into an equivalent Eulerian 
density $\f(\x,t)$.
The discretized version of this operation is called \emph{force spreading}. 
We express the force spreading operation by $\f = \cS[\X] \, \F$, in which $\cS[\X]$ is the \emph{force-spreading operator}.
Eq.~\eqref{eqn-u-interpolation} determines the physical velocity of each Lagrangian material point from the Eulerian velocity field, so that the immersed structure moves according to the local value of the velocity field $\u(\x,t)$.
This \emph{velocity interpolation} operation is expressed as $\D{\X}{t} = \cJ[\X] \, \u$, in which $\cJ[\X]$ is the \emph{velocity-interpolation operator}.
Notice that $\cS$ and $\cJ$ are adjoint operators, $\cS = \cJ^{*}$ \cite{Peskin02}.

\subsection{Discrete equations of motion}

We consider only linear solver performance in this work using linear systems of equations that arise from an energy-stable semi-implicit discretization of the IB equations.
We use a spatial discretization that is similar to one used in earlier work \cite{RDGuy15-gmgiib}, which is briefly described in \ref{sec-spatial-discretization}.
To discretize these equations in time, let $\dt$ be the time step size, and let $n$ be the time step number.
In each time step, we simultaneously solve for the updated Eulerian velocity $\u^{n+1}$ and pressure $p^{n+1}$ at time $t^{n+1} = (n+1) \dt$ along with the structural configuration $\X^{n+1}$.
To simplify notation, we use $\cS_h^{n} \equiv \cS_h[\X^{n}]$ to indicate the spreading operator corresponding to structural configuration $\X^{n}$ along with analogous notation for the interpolation operator $\cJ_h$.
The time-stepping scheme reads
\begin{align}
	\rho \left(\frac{\u^{n+1} - \u^n}{\dt} + [\u \cdot \grad_h \u]^{(n+\half)}\right) &= -\grad_h p^{n+1} + \mu \V{\grad}^2_h \u^{n+1} + \cS_h^{n} \F^{n+1}, \label{eqn-imp-momentum} \\
	\grad_h \cdot \u^{n+1} &= \V 0, \label{eqn-imp-continuity} \\
	\frac{\X^{n+1} - \X^n}{\dt} &= \cJ_h^{n} \u^{n+1}, \label{eqn-imp-X-new} \\
	\F^{n+1} &= \cF_h\left[\X^{n+1}\right]. \label{eqn-imp-F}
\end{align}
Except for the nonlinear convection term, this scheme uses a combination of forward and backward Euler time stepping.
In all tests reported herein, we omit the convection term, because including it does not affect linear solver performance.
In applications, however, we often use a version of the PPM method \cite{ColellaWoodward84, RiderEtAl07, BEGriffith09-efficient} along with Adams-Bashforth to approximate the midstep value of $\u \cdot \grad \u$ via
\begin{equation}
	[\u \cdot \grad_h \u]^{(n+\half)} = \frac{3}{2} \u^{n} \cdot \grad_h \u^{n} - \frac{1}{2} \u^{n-1} \cdot \grad_h \u^{n-1}.
\end{equation}

Under reasonable assumptions on the form of the discretized force operator $\cF_h$, the only stability restriction associated with this semi-implicit time stepping scheme is related to our explicit treatment of the convective term in the momentum equation \cite{Newren07}.
A fully implicit version of this scheme would replace $\cS_h^{n} = \cS_h[\X^{n}]$ by $\cS_h^{n+1} = \cS_h[\X^{n+1}]$, and likewise for $\cJ_h$.
Such schemes do not appear to offer benefits in terms of energy stability or order of accuracy~\cite{Newren07}, but they do require the solution of a more complex system of nonlinear equations.
The semi-implicit formulation used here can be seen as a method that linearizes the geometrical nonlinearities associated with the coupling operators $\cS_h$ and $\cJ_h$ without sacrificing energy stability or formal order of accuracy.

Depending upon the functional form of the discrete force operator $\cF_h$, Eqs.~\eqref{eqn-imp-momentum}--\eqref{eqn-imp-F} can be linear or nonlinear. 
Because we consider only linear solver performance in this work, we use only linear force functionals of the form $\cF_h[\X] = \cK_h \X$, in which $\cK_h$ is the stiffness matrix of the elasticity model, and we solve Eqns.~\eqref{eqn-imp-momentum}--\eqref{eqn-imp-F} by a preconditioned Krylov method.
This requires the solution of linear systems of the form
\begin{equation}
	\left(\begin{array}{ccc}
		\cA & \cG  & -\cS_h^{n} \cK_h  \\
		-\cD & \V 0  & \V 0 \\
		-\cJ_h^{n} & \V 0 & \frac{1}{\dt}\V{I}
	\end{array}\right)
	\left(\begin{array}{c}
		\u^{n+1}\\
		p^{n+1} \\
		\X^{n+1}
	\end{array}\right) =
	\left(\begin{array}{c}
		\V{g} \\
		0 \\
	 	\frac{1}{\dt} \X^{n}
	\end{array}\right), \label{eqn-system-impib}
\end{equation}
in which $\cA = \frac{\rho}{\dt}\I - \mu \V{\grad}^2_h$, $\cG = \grad_h$, and $\cD = \grad_h \cdot \mbox{}$ are block Eulerian operators, and $\V{g}$ contains contributions from previous time steps and explicitly evaluated terms from the current time step.
To develop a system of equations amenable to solution via geometric multigrid methods, we use Eq.~\eqref{eqn-imp-X-new} to eliminate $\X^{n+1}$ from the block system and obtain
\begin{equation}
	\underbrace{\left(\begin{array}{cc}
		\cAIB   &  \cG   \\
		-\cD      & \V 0
	\end{array}\right)}_{\cLIB}
	\left(\begin{array}{c}
		\u^{n+1}\\
		p^{n+1}
	\end{array}\right) =
	\left(\begin{array}{c}
		\V{g} + \cS_h^{n} \cK_h \X^n \\
		\V 0
	\end{array}\right), \label{eqn-matrix-impib}
\end{equation}
in which $\cAIB = \cA - \dt \cS_h^{n} \cK_h \cJ_h^{n}$ is the modified momentum operator that includes the projection of $\cK_h$, the linear 
Lagrangian elasticity operator, onto the Eulerian frame.
We refer to $\cLIB$ as the \emph{Stokes-IB operator}.


\section{Multigrid}
\label{sec-MG}

\subsection{Basic multigrid algorithm}

An effective preconditioner is needed to solve the Stokes-IB system~\eqref{eqn-matrix-impib} efficiently using a Krylov method.
Here, we briefly discuss the key ingredients of a geometric multigrid (GMG) preconditioner for the Stokes-IB system.
Detailed descriptions of the multigrid method are available~\cite{BriggsHensonMcCormick00, UTrottenberg01}, and the development of GMG methods for the IB method is also discussed in previous work~\cite{RDGuy15-gmgiib, RDGuy12}.

We construct a hierarchy of uniform Cartesian discretizations of the spatial domain $\Omega$.
Let $\Omega^\ell$ indicate a particular discretization with grid spacing $h^\ell$, in which $\ell = 0, 1, \ldots, \ellmax$ indicates the \emph{level} of the discretization, with $\ell = 0$ denoting the coarsest level in the hierarchy and $\ell = \ellmax$ denoting the finest level.
The grid spacings on adjacent levels $\ell$ and $\ell - 1$ are related by an integer refinement ratio $r_{\text{ref}}$, so that $h^\ell = \frac{h^{\ell-1}}{r_{\text{ref}}}$.
Here, we only consider $r_{\text{ref}} = 2$.

\begin{algorithm}[t]
	\caption{recursive V-cycle multigrid}
    \begin{algorithmic}[1]
	\Procedure{\bf $\w^\ell \leftarrow \textbf{MG}$}{$\w^\ell, \b^\ell, \Omega^\ell, \nu_1, \nu_2$}
	\If {$\Omega^\ell$ is the coarsest grid}
		\State $\w^\ell \leftarrow (\cLIB^\ell)^{-1} \b^\ell$  \Comment{solve the coarse grid equation}
	\Else
		\State $\w^\ell \leftarrow \textbf{smooth}(\w^\ell, \b^\ell, \nu_1)$  \Comment{apply $\nu_1$ pre-smoothing sweeps}
		\State $\r^\ell \leftarrow \b^\ell - \cLIB^\ell \w^\ell$  \Comment{compute the residual on the present level}
		\State $\r^{\ell-1} \leftarrow \cR_\ell^{\ell-1} \r^\ell$  \Comment{restrict the residual to the next coarser level}
		\State $\e^{\ell-1} \leftarrow \textbf{MG}(\V 0, \r^{\ell-1}, \Omega^{\ell-1}, \nu_1, \nu_2)$ \Comment{recursively call \textbf{MG}}
		\State $\e^\ell \leftarrow \cP_{\ell-1}^\ell \e^{\ell-1}$  \Comment{prolong the error from the next coarser level}
		\State $\w^\ell \leftarrow \w^\ell + \e^\ell$  \Comment{correct the solution on the present level}
		\State $\w^\ell \leftarrow \textbf{smooth}(\w^\ell, \b^\ell, \nu_2)$  \Comment{apply $\nu_2$ post-smoothing sweeps}
	\EndIf
    \EndProcedure
\end{algorithmic}
\label{alg-MG-V}
\end{algorithm}

The basic V-cycle multigrid algorithm used in this work is shown in Algorithm~\ref{alg-MG-V}.
It aims to solve the discretized equations on $\Omega^{\ell}$ by combining simple approximate solvers on level $\ell$ with coarse-grid corrections recursively computed on levels $\ell-1, \ell-2, \ldots$.
Specifically, on each grid level $\ell > 0$, Algorithm~\ref{alg-MG-V} uses a \emph{smoother} to eliminate the high-frequency components of the error.
The remaining low-frequency errors are meant to be eliminated by coarse-grid corrections.
(Although not shown here, alternative multigrid algorithms, such as F-~or W-cycles, seem to offer little to no benefit for the Stokes-IB equations with our present smoothers.)

Grid levels are connected by a \emph{restriction operator} $\cR_\ell^{\ell-1}$ that coarsens solution data from a finer level $\ell$ to a coarser level $\ell-1$, and a \emph{prolongation operator} $\cP_{\ell-1}^\ell$ that interpolates solution data from a coarser level $\ell-1$ to a finer level $\ell$.
These operators are also used to define the Stokes-IB system on coarser grid levels.
Omitting the dependence on the time step number $n$, the block linear system on level $\ell$ is
\begin{equation}
	\left(\begin{array}{cc}
		\cA^\ell -\dt \left[\cS_h \cK_h \cJ_h\right]^\ell	& \cG^\ell \\
		-\cD^\ell												& \V 0
	\end{array}\right)
	\left(\begin{array}{c}
		\u^\ell \\
		p^\ell
	\end{array}\right) =
	\left(\begin{array}{c}
		\V{g}^\ell \\
		\V 0
	\end{array}\right). \label{eqn-matrix-h-impib}
\end{equation}
The coarse-grid Eulerian elasticity operator is defined for $\ell < \ellmax$ via a Galerkin projection,
\begin{equation}
	\left[\cS_h \cK_h \cJ_h\right]^\ell = \left(\cR_{\u}\right)_{\ell+1}^\ell \left[\cS_h \cK_h \cJ_h\right]^{\ell+1} \left(\cP_{\u}\right)_\ell^{\ell+1},
\end{equation}
in which $\cR_{\u}$ and $\cP_{\u}$ are restriction and prolongation operators for velocity-like degrees of freedom only.
The coarse-grid operators $\cA^\ell$, $\cG^\ell$, and $\cD^\ell$ are constructed by rediscretization.
We denote the block system~\eqref{eqn-matrix-h-impib} by
\begin{equation}
	\cLIB^\ell \w^\ell = \b^\ell. \label{eqn-mg-h}
\end{equation}
Given an approximate solution $\widetilde{\w}^\ell$, the corresponding error equation is
\begin{equation}
	\cLIB^\ell \e^\ell = \r^\ell, \label{eqn-mg-error-h}
\end{equation}
in which $\e^\ell = \w^\ell - \widetilde{\w}^\ell$ is the error and $\r^\ell = \b^\ell - \cLIB^\ell \widetilde{\w}^\ell$ is the residual.

\subsection{Prolongation and restriction}

The velocity prolongation operator $\cP_{\u}$ is based on lowest-order Raviart-Thomas interpolation on quadrilaterals \cite{DNArnold05}, which uses piecewise-linear interpolation in the normal direction to the cell face along with piecewise-constant interpolation in the tangential direction.
The velocity restriction operator is taken to be the adjoint of the prolongation operator, $\cR_{\u} = \cP_{\u}^{*}$.
Consequently, the coarse-grid versions of the Eulerian elasticity operator $\cS_h \cK_h \cJ_h$ retain the symmetry of the fine-grid operator.
We also use linear interpolation to prolong pressure values from coarse to fine grid levels, and we use simple averaging to restrict pressures from fine to coarse levels.
In this case $\cR_p \neq \cP_p^{*}$, but preliminary tests suggest that this has essentially no effect on solver convergence.

\subsection{Smoothers}

We consider two classes of smoothers: Vanka-like algorithms based on Schwarz domain decomposition methods, and algorithms based on an approximate block factorization.
In both cases, the basic algorithm may not act as a smoother (i.e.~it may fail to damp some high-frequency error components).
To provide enhanced smoothing without introducing a damping parameter, we use a fixed number of flexible GMRES (FGMRES) \cite{Saad93} iterations preconditioned by the basic algorithm.

\subsubsection{Schwarz smoothing}

Additive and multiplicative Schwarz are domain decomposition methods~\cite{SmithBjorstadGropp96} that solve restricted versions of the linear system on overlapping subdomains.
The key difference between additive and multiplicative domain decomposition methods is that, in an additive method, the subdomain-restricted equations are solved independently, whereas in a multiplicative method, the most recently computed solution values from all subdomains are used in computing the residual for each subdomain solve.
Additive and multiplicative Schwarz are thereby generalizations of the classical Jacobi and Gauss-Seidel methods.
Multiplicative algorithms are fundamentally sequential, making them difficult or impossible to parallelize.
Additive algorithms, by contrast, are readily parallelized.

  \begin{algorithm}[t]
  \caption{restricted additive Schwarz (RAS)}\label{alg-RAS}
  \begin{algorithmic}[1]
    \Procedure{$\w \leftarrow \textbf{RAS}$}{$\w,\b$}
      \State partition $\Omega$ into $N$ overlapping subdomains $\Omega_i$ with index sets $G_i$ and restriction matrices $R_i$, and construct non-overlapping index sets $\tilde{G}_i \subset G_i$ and corresponding restriction matrices $\tilde{R}_i$
      \State $\r \gets \b - \cL \w$ \Comment{form the residual for the \emph{initial} value of $\w$}
      \For{partition $i = 1 \dots N$}
        \State $\cL_i \gets R_i \cL R_i^{T}$ \Comment{construct the subdomain operator $\cL_i$}
        \State $\b_i \gets R_i \r$ \Comment{form the local right-hand side $\b_i$}
        \State $\w \gets \w + \tilde{R}_i^{T} (\cL_i)^{-1}\b_i$ \Comment{perform a local solve and update $\w$ in $\tilde{G}_i$}
      \EndFor
    \EndProcedure
  \end{algorithmic}
  \end{algorithm}

  \begin{algorithm}[t]
  \caption{restricted multiplicative Schwarz (RMS)}\label{alg-RMS}
  \begin{algorithmic}[1]
    \Procedure{$\w \leftarrow \textbf{RMS}$}{$\w,\b$}
      \State partition $\Omega$ into $N$ overlapping subdomains $\Omega_i$ with index sets $G_i$ and restriction matrices $R_i$, and construct non-overlapping index sets $\tilde{G}_i \subset G_i$ and corresponding restriction matrices $\tilde{R}_i$
      \For{partition $i = 1 \dots N$}
        \State $\r \gets \b - \cL \w$ \Comment{form the residual for the \emph{updated} value of $\w$}
        \State $\cL_i \gets R_i \cL R_i^{T}$ \Comment{construct the subdomain operator $\cL_i$}
        \State $\b_i \gets R_i \r$ \Comment{form the local right-hand side $\b_i$}
        \State $\w \gets \w + \tilde{R}_i^{T} (\cL_i)^{-1}\b_i$ \Comment{perform a local solve and update $\w$ in $\tilde{G}_i$}
      \EndFor
    \EndProcedure
  \end{algorithmic}
  \end{algorithm}

In this work, we consider two Schwarz-like algorithms.
One is \emph{restricted additive Schwarz (RAS)} \cite{XCCai99, EEfstathiou03-RAS} (Algorithm~\ref{alg-RAS}) with overlapping subdomains that correspond to all of the degrees of freedom associated with contiguous, rectangular boxes of grid cells. 
The other is a multiplicative version of the RAS algorithm, which we refer to as \emph{restricted multiplicative Schwarz (RMS)} (Algorithm~\ref{alg-RMS}).
Notice that the only difference between the two algorithms is the manner in which the residual is computed for each of the subdomain solves.
Our RMS smoother is similar to the ``big-box'' Vanka algorithm described by Guy et al.~\cite{RDGuy15-gmgiib}.
We consider the effect of different subdomain sizes and overlap widths on the performance of both algorithms.


\subsubsection{Schur complement (SC) smoothing}

An alternative smoothing approach is to construct an approximate block factorization of $\left(\cLIB^\ell\right)^{-1}$ that separates the block system into velocity and pressure subdomain operators.
Omitting the superscript ``$\ell$'' for notational clarity, we first note that $\cLIB^{-1}$ can be written as
 \begin{equation}
\cLIB^{-1} = \left( \begin{array}{cc}
\I & \; - \cAIB^{-1} \cG \\
\V 0 & \I
\end{array} \right)  \left(\begin{array}{cc}
\cAIB^{-1} & \V 0\\
\V 0 & \cM^{-1}
\end{array}\right)  \left(\begin{array}{cc}
\I & \V 0\\
\cD \cAIB^{-1} & \I
\end{array}\right),\label{eqn-fact-Stokes-IB}
\end{equation}
in which $\cM = \cD \cAIB^{-1} \cG$ is the \emph{Schur complement} of the block system of equations in~\eqref{eqn-matrix-impib}.
The cost of constructing the full factorization of $\cLIB^{-1}$ is prohibitive, but in practice, $\cLIB^{-1}$ rarely needs to be formed explicitly.
Here, we only need to be able to apply this operator to known right-hand side vectors.
Moreover, the application of the exact inverse operator can be unnecessary in a smoother, since we only wish to eliminate high-frequency error modes and not low-frequency modes.
Consequently, in a multigrid algorithm, operators that approximate the action of $\cAIB^{-1}$ and $\cM^{-1}$ can suffice.

There are many choices for approximating $\cAIB^{-1}$ and $\cM^{-1}$, but we find that simple approximate solvers for $\cAIB$ and $\cM$ lead to an effective smoothing algorithm.
We approximate the action of $\cAIB^{-1}$ by $\ctAIB^{-1}$, which uses a fixed number of Chebyshev iterations preconditioned by Gauss-Seidel applied to $\cAIB$, and we approximate the action of $\cM^{-1}$ by $\ctM^{-1}$, which uses a fixed number of Chebyshev iterations for the operator $\left(\cD \ctAIB^{-1} \cG\right)$ preconditioned by Gauss-Seidel applied to a sparse approximate Schur complement, $\chM = \cD \; \left(\text{diag}(\cAIB)\right)^{-1} \; \cG$.
Notice that the sparse approximate Schur complement $\chM$ takes the form of an inhomogeneous discrete Poisson operator.
We typically use two Chebyshev iterations for both $\ctAIB^{-1}$ and $\ctM^{-1}$.
With these operators so defined, we specify the action of $\ctLIB^{-1} \approx \cLIB^{-1}$ by
\begin{equation}
	\ctLIB^{-1} =
	\left( \begin{array}{cc}
		\I & \; - \ctAIB^{-1} \cG \\
		\V 0 & \I
	\end{array} \right)
	\left(\begin{array}{cc}
		\ctAIB^{-1} & \V 0\\
		\V 0 & \ctM^{-1}
	\end{array}\right)
	\left(\begin{array}{cc}
		\I & \V 0\\
		\cD \ctAIB^{-1} & \I
	\end{array}\right).\label{eqn-approx-fact-Stokes-IB}
\end{equation}
Notice that this smoother involves only point-relaxation.
By contrast, in two spatial dimensions, the Vanka-like smoothers relax over Cartesian boxes of size $n_x \times n_y$, which requires solving linear systems with $O(n_x n_y)$ variables.
Consequently, the typical cost of an application of the SC smoother is substantially less than the cost of one application of the Vanka-like smoother.

\section{Software implementation}

The solvers for this semi-implicit IB method are implemented in the open-source IBAMR library~\cite{IBAMR-web-page}.
We use PETSc \cite{petsc-efficient, petsc-user-ref, petsc-web-page} to provide Krylov solvers 
and Schur complement-based smoothers.
A custom implementation of the Schwarz algorithms is provided by IBAMR.
IBAMR also relies on SAMRAI \cite{HornungKohn02, samrai-web-page} for Cartesian grid management, and for inter-level data transfer operations.

\section{Results}

Our tests explore the linear solver performance for both thick and thin elastic structures at a range for material stiffnesses and flow conditions.
We consider immersed structures that are modeled using systems of elastic fibers.
For simplicity, we assume that the leading Lagrangian coordinate $s_1$ varies along the direction of each fiber, and that the remaining curvilinear coordinates serve to label a particular fiber.
The unit fiber tangent vector is $\V{\tau} = \D{\X}{s_1}/ \left\|\D{\X}{s_1}\right\|$, and the tension $T$ in each fiber is a function of the fiber strain $\left\|\D{\X}{s_1}\right\|$.
It can be shown that the Lagrangian fiber force functional takes the form \cite{Peskin02}
\begin{equation}
 \F = \cF[\X] = \D{(T \V{\tau})}{s_1}. \label{eqn-F-functional}
\end{equation}
To explore the performance of the linear solvers, it is convenient to use fibers with zero resting lengths that resist only extension, so that the fiber tension is $T = \alpha \left\|\D{\X}{s_1}\right\|$, in which $\alpha$ is the fiber stiffness.
The resulting Lagrangian force density is
\begin{equation}
 \F = \alpha \DD{\X}{s_1}, \label{eqn-F-lin-functional}
\end{equation}
which is a linear functional.
We use a simple second-order finite difference approximation to this functional, as described in \ref{sec-lagrangian-discretization}.

In all of our tests, the physical domain is $\Omega = [0, 1]^2$, and we use regularized lid-driven cavity flow boundary conditions, for which the velocity is set to zero along $\partial\Omega$ except along the top wall, where we prescribe $u(x,1) = (1 - \cos (2 \pi x))/2$ and $v = 0$.
In these tests, the physical boundary conditions set a characteristic flow speed that determines the dynamic timescale of interest, independent of the elastic timescales of the immersed structure.
The domain size sets the characteristic lengthscale as $L = 1$, and lid-driven cavity flow conditions set the characteristic flow speed as $U = 1$, so that $Re = \frac{\rho U L}{\mu} = \frac{\rho}{\mu}$.
We set $\rho = 1$ for nonzero Reynolds number cases, yielding $\mu = Re^{-1}$.

We only consider square computational domains, and in our multigrid algorithm, we always use an $8 \times 8$ coarse grid along with sufficiently many finer levels to reach the targeted Eulerian grid spacing.
A direct solver is used on the coarsest grid level.
As mentioned previously, we only consider $r_{\text{ref}} = 2$.

\subsection{Thick elastic shell}
\label{sec-linear-shell}

\begin{figure}[t]
\centering
\includegraphics[width=0.625\textwidth]{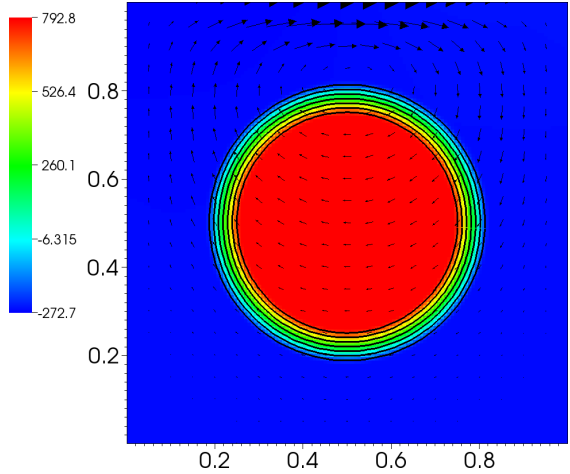}
\caption{
	Representative results from the thick elastic test problem of Sec.~\ref{sec-linear-shell}.
	Here, we plot the velocity and pressure field along with the structure configuration for $\mu=0$ and a relative stiffness of $\gamma=5$.
	The Eulerian grid spacing is $h = \frac{1}{128}$.
	For clarity, only half of the Lagrangian fibers are shown.
	Results for different values of $\mu$ and $\gamma$ are similar.
}
\label{fig-steady_shell}
\end{figure}

We first consider a thick circular annulus described using Lagrangian curvilinear coordinates $(s_1,s_2) \in U = [0, 2 \pi) \times [0, w]$, in which $w = 1/16$ is the thickness of the annulus.
The structure is initially placed at the center of the domain in the configuration
\begin{equation}
 	\X(s_1, s_2) = \left(\xc + (r + s_2)\cos(s_1), \yc + (r + s_2)\sin(s_1) \right). \label{eqn-annulus}
\end{equation}
We choose the center to be $\x_\text{c} = (\xc, \yc) = (0.5, 0.5)$ and the inner radius to be $r = 1/4$.
This configuration has been used as a standard test case in the IB literature~\cite{RDGuy15-gmgiib, DBoffi08, BEGriffith05-ib_accuracy, BEGriffith07-ibamr_paper, BEGriffith12-ib_volume_conservation}.
The Eulerian domain is discretized using a uniform $N \times N$ grid, so that the Cartesian grid spacing is $h = \dx_1 = \dx_2 = 1/N$.
For the Lagrangian domain, we use $M_1 = \frac{19}{8} N$ points in the $s_1$ direction and $M_2 = \frac{3}{32} N +1$ points in the $s_2$ direction, which yields a physical spacing between the Lagrangian nodes approximately equal to $\frac{2}{3}h$.
In these tests, we always use $\dt = 0.32 \dx$.
This implies that the time step size satisfies a mild advective CFL-type condition under grid refinement.

We characterize the stiffness $\alpha$ in terms of a stiffness ratio $\gamma$ via
\begin{equation}
	\alpha = \gamma \frac{3.93}{0.005}.
\end{equation}
In prior work \cite{RDGuy15-gmgiib}, $\alpha = \frac{3.93}{0.005}$ was shown to be approximately the largest stiffness for this problem for which the explicit solver is stable for a time step size of $\dt = 0.005$ and a grid spacing of $\dx = \frac{1}{64}$ in Stokes flow conditions.
This largest stiffness is relatively insensitive to grid refinement \cite{RDGuy15-gmgiib}.
Thus, $\gamma$ roughly characterizes the ratio of the stiffness $\alpha$ and the maximum stiffness that can be used by an explicit time stepping scheme, $\alpha_\text{explicit}$.
For non-zero Reynolds numbers, $\alpha_\text{explicit}$ also depends on the fluid viscosity $\mu$ in a manner that we do not explore in this work.
Thus, $\gamma$ only \emph{approximately} characterizes the ratio of the elastic stiffness to the largest stiffness that can be used by an explicit time stepping scheme.
We consider a range of values of $\gamma$, from relatively soft ($\gamma = 5$) to very stiff ($\gamma = 500$).
Representative results are shown in Fig.~\ref{fig-steady_shell}.

\subsubsection{Schwarz smoothers}

\begin{figure}[t]
\centering
\includegraphics[width=0.6\textwidth]{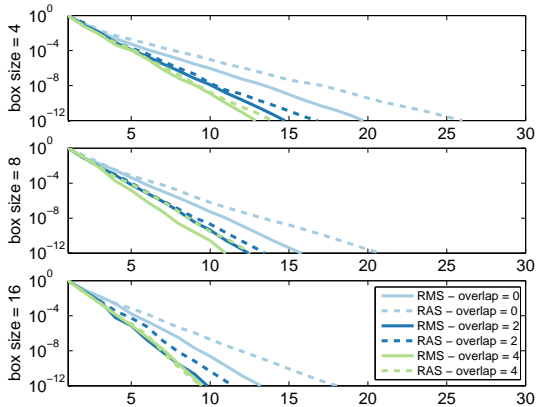}
\caption{
	Multiplicative (solid line) and additive (dashed line) smoother performance under Stokes flow conditions for the thick elastic shell (Sec.~\ref{sec-linear-shell}) with $N = 128$ and $\gamma = 500$, and for different subdomain sizes and overlap widths.
	For this problem, solver performance is relatively insensitive to subdomain size, overlap width, or smoother algorithm.
}
\label{fig-5.1.1-1}
\end{figure}

\begin{figure}
\centering
\footnotesize
  \begin{tabular}{>{\centering\arraybackslash}m{32pt} >{\centering\arraybackslash}m{0.8\textwidth} }
  (a) \vspace{-30pt} & \\ $\mu = 1.0$   & {\includegraphics[width=0.6\textwidth]{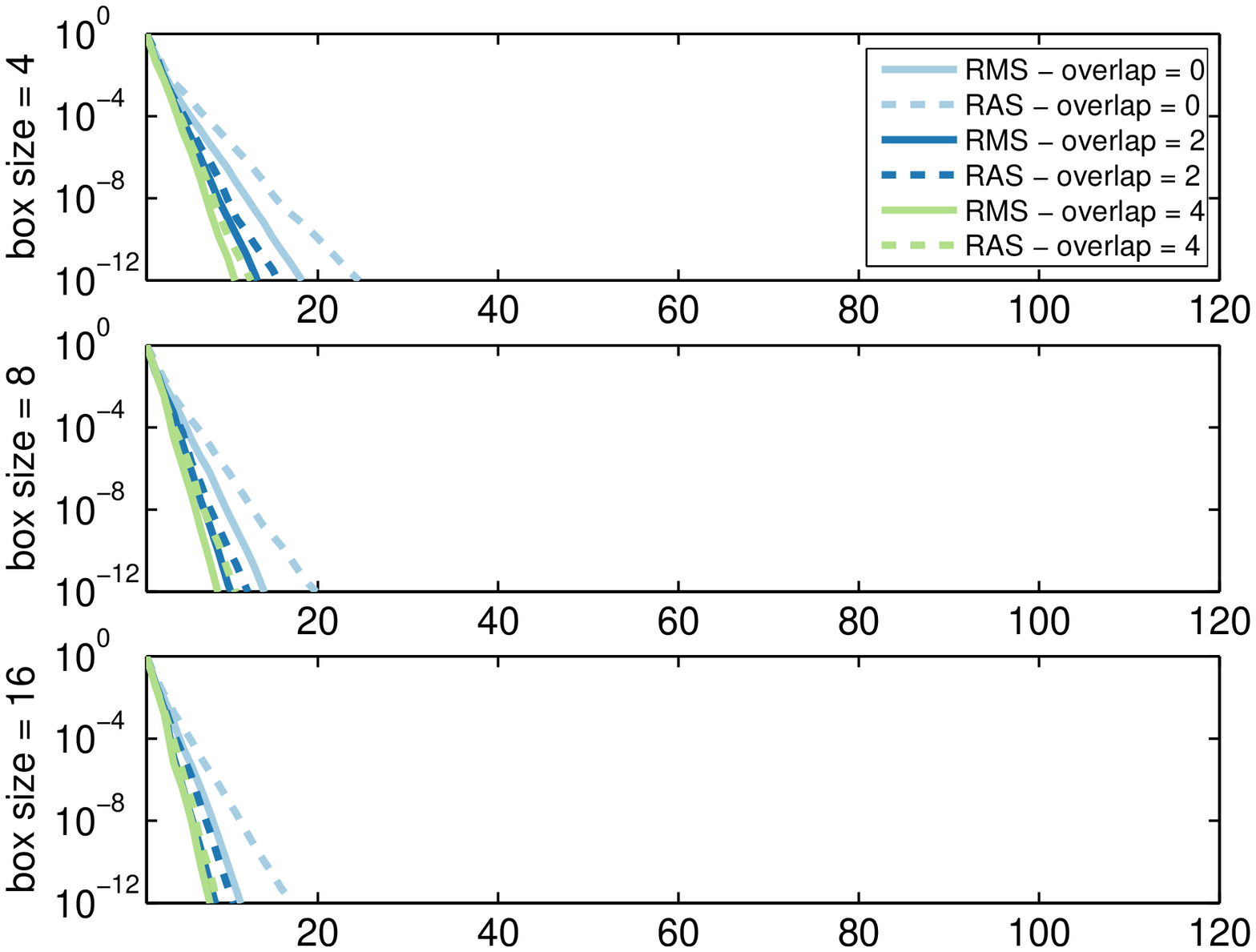}} \vspace{-1.1\baselineskip} \\
  (b) \vspace{-30pt} & \\ $\mu = 0.1$   & {\includegraphics[width=0.6\textwidth]{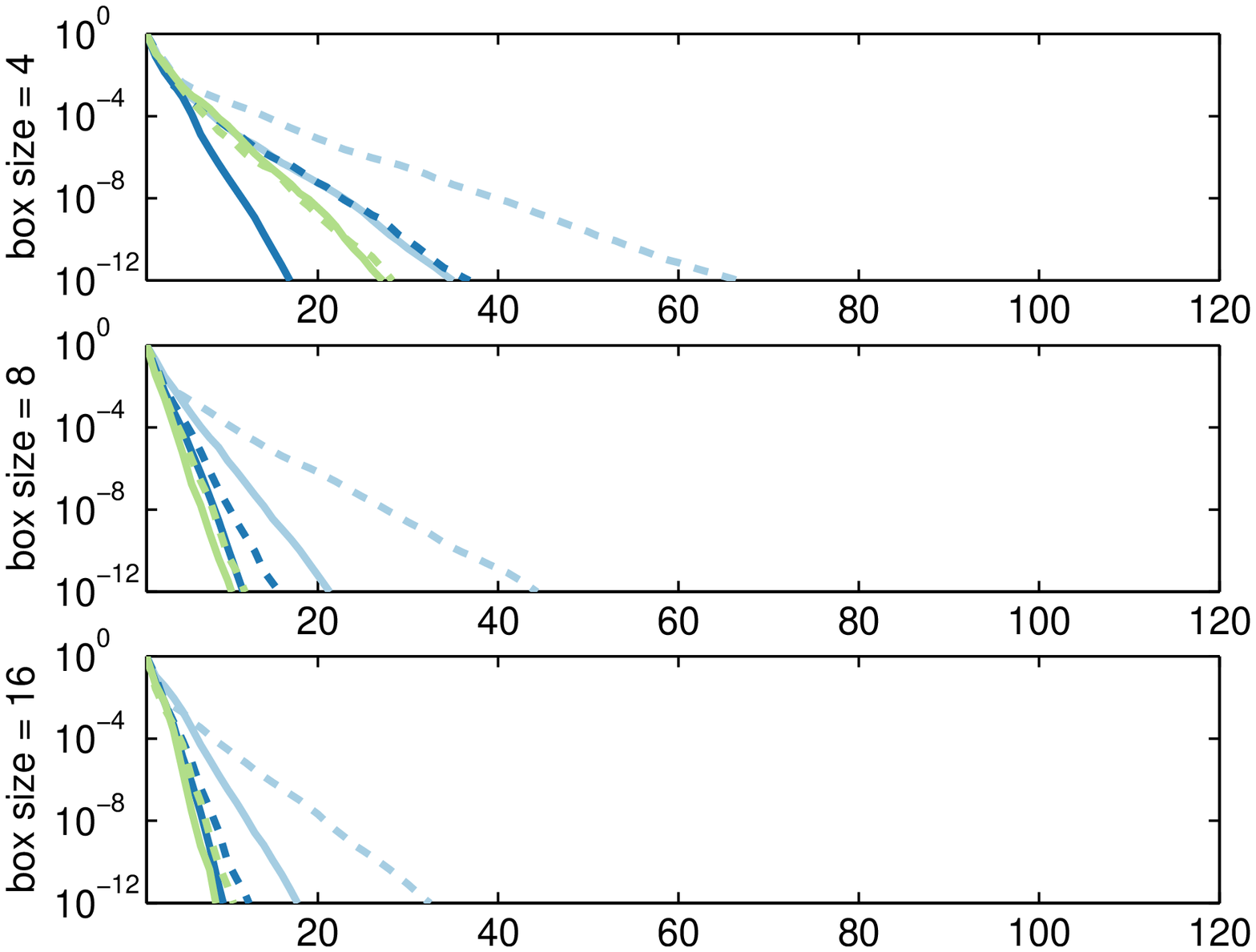}} \vspace{-1.1\baselineskip} \\
  (c) \vspace{-30pt} & \\ $\mu = 0.01$  & {\includegraphics[width=0.6\textwidth]{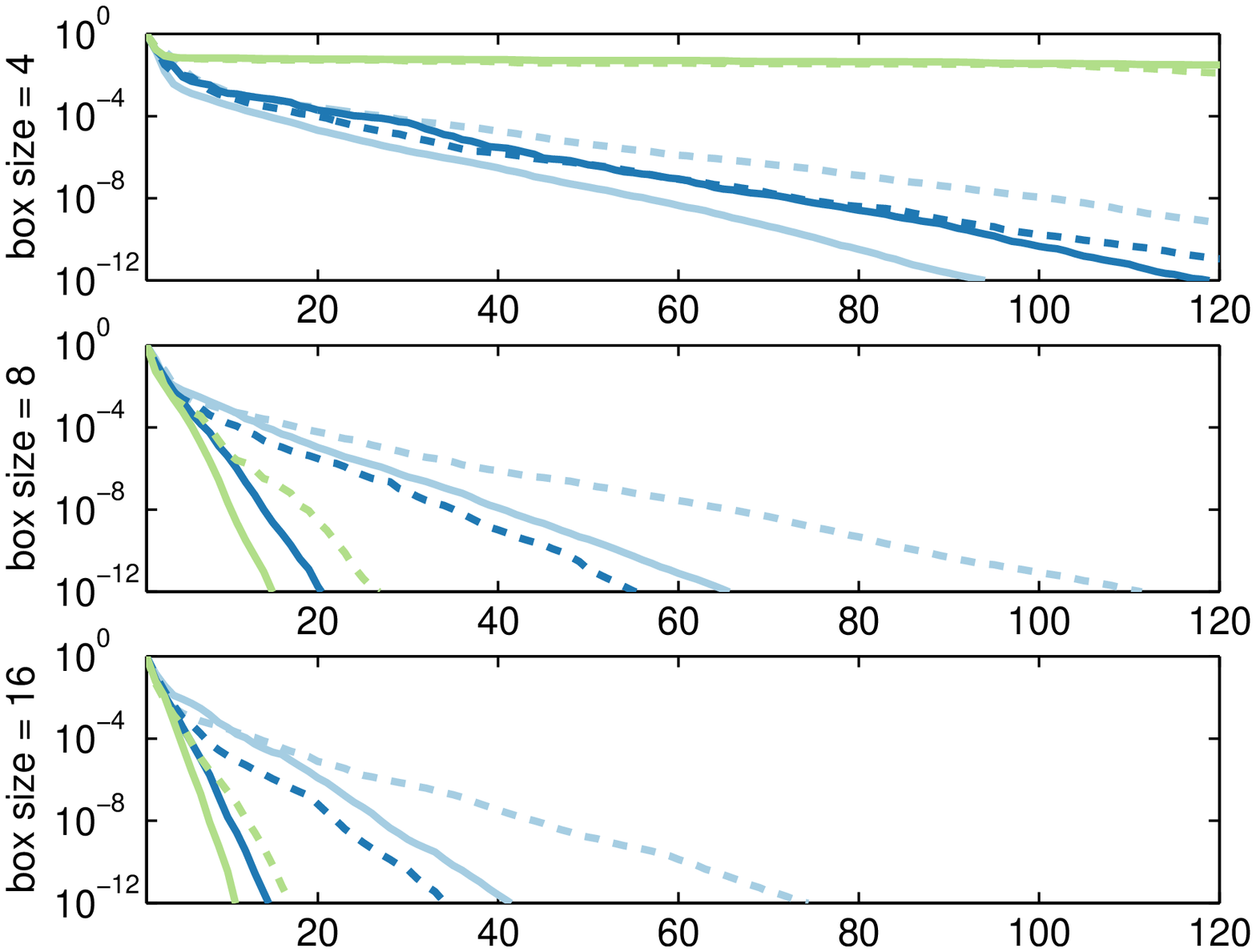}}
  \end{tabular}
\caption{
	Similar to Fig.~\ref{fig-5.1.1-1}, but here showing the effect of decreasing viscosity, with (a)~$\mu = 1$, (b)~$\mu = 0.1$, and (c)~$\mu = 0.01$.
	Performance clearly degrades with decreasing viscosity.
	In the most challenging cases, the multiplicative smoother outperforms the additive smoother by a wide margin.
}
\label{fig-5.1.1-2}
\end{figure}

\begin{figure}
\centering
\footnotesize
  \begin{tabular}{>{\centering\arraybackslash}m{32pt} >{\centering\arraybackslash}m{0.8\textwidth} }
  (a) \vspace{-30pt} & \\ RMS & {\includegraphics[width=0.75\textwidth]{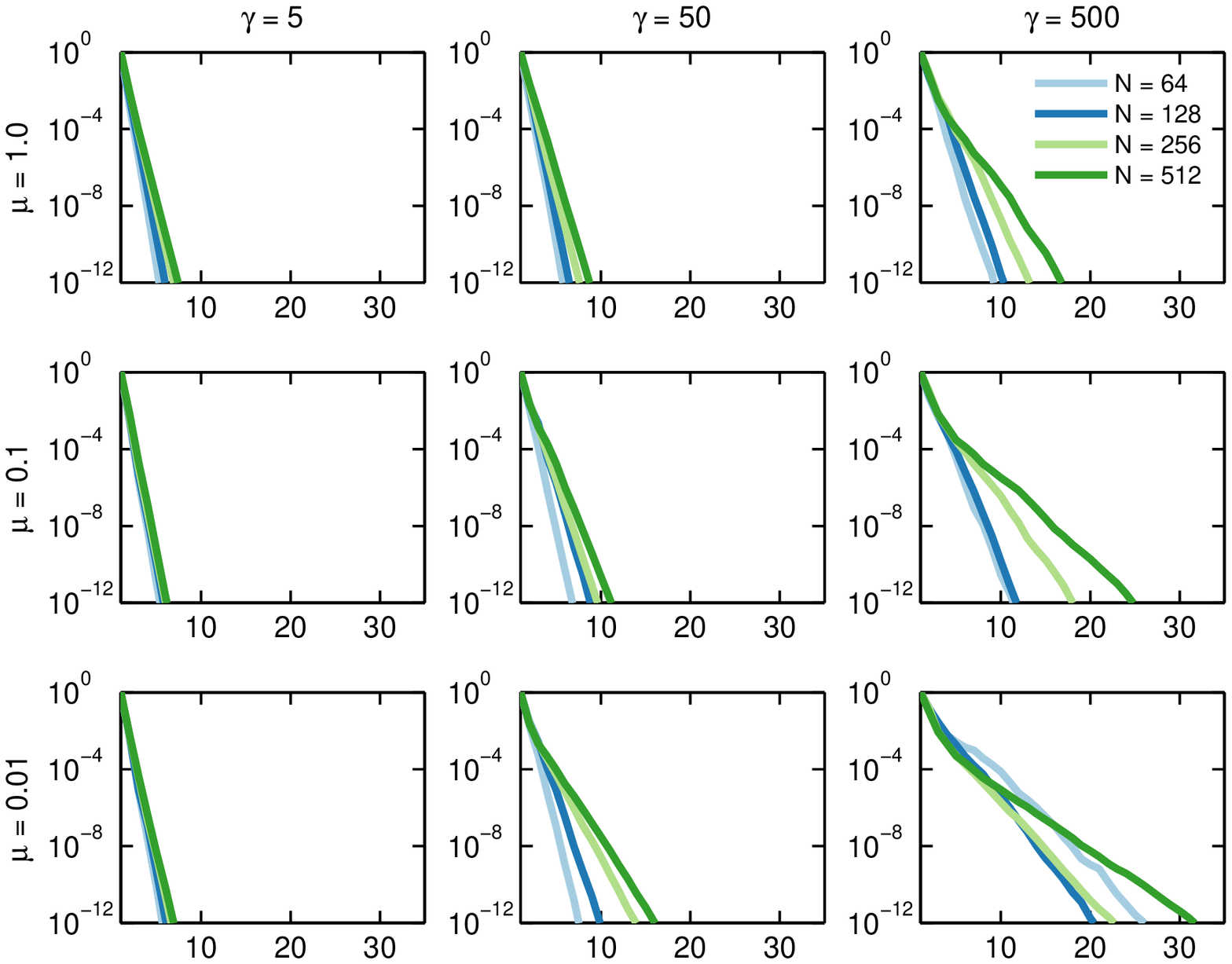}} \vspace{-1.1\baselineskip} \\
  (b) \vspace{-30pt} & \\ RAS & {\includegraphics[width=0.75\textwidth]{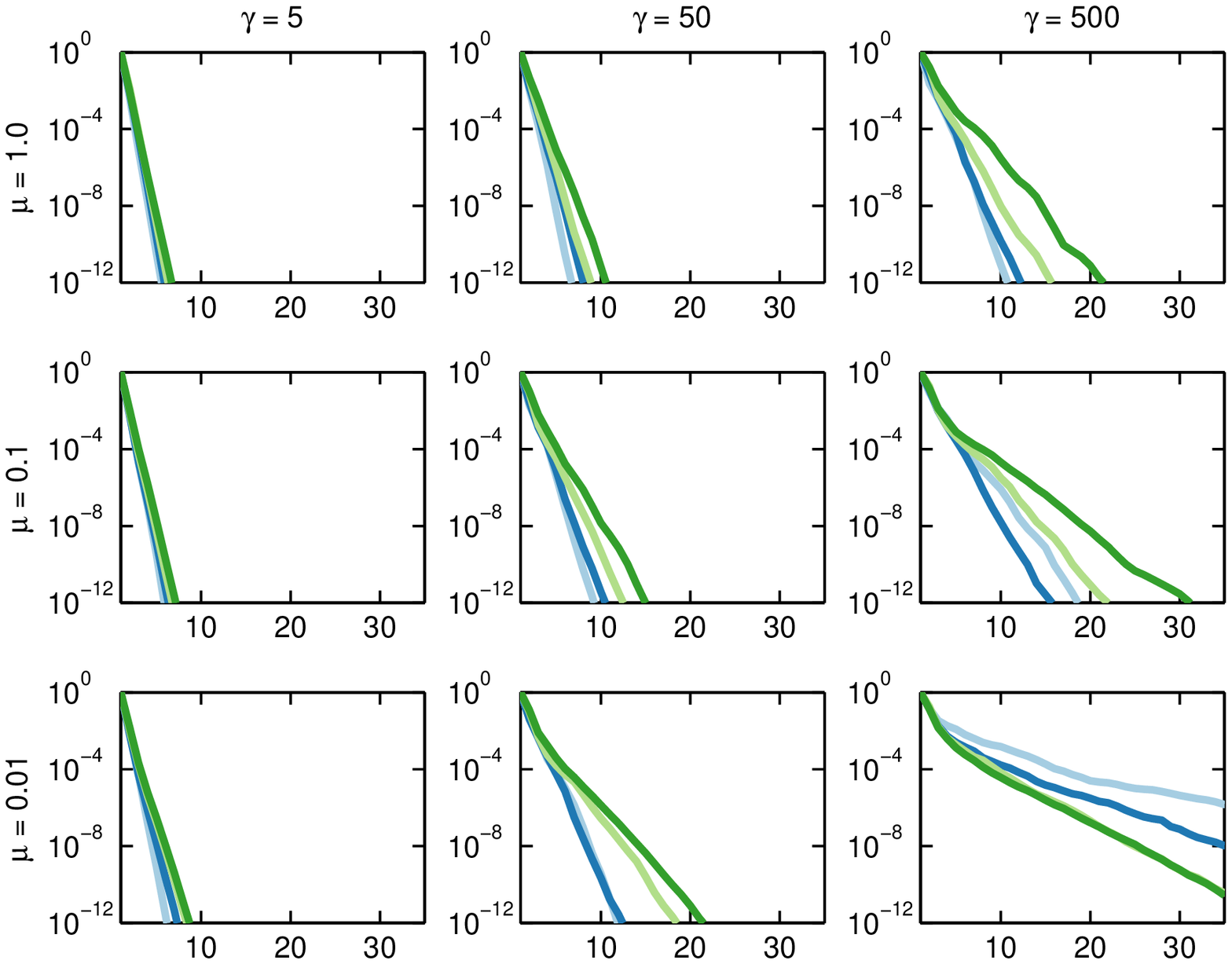}}
  \end{tabular}
\caption{
	Performance of multigrid using (a) multiplicative (RMS) and (b) additive (RAS) smoothers under grid refinement for a range of relative stiffnesses ($\gamma$) and viscosities ($\mu$), using subdomains of size $8 \times 8$ and an overlap width of $2$.
	In most cases, using the additive smoother results in only a modest increase in iterations compared to the multiplicative algorithm.
}
\label{fig-5.1.1-3}
\end{figure}

\begin{figure}
\centering
\footnotesize
  \begin{tabular}{>{\centering\arraybackslash}m{32pt} >{\centering\arraybackslash}m{0.8\textwidth} }
  (a) \vspace{-30pt} & \\ RMS & {\includegraphics[width=0.75\textwidth]{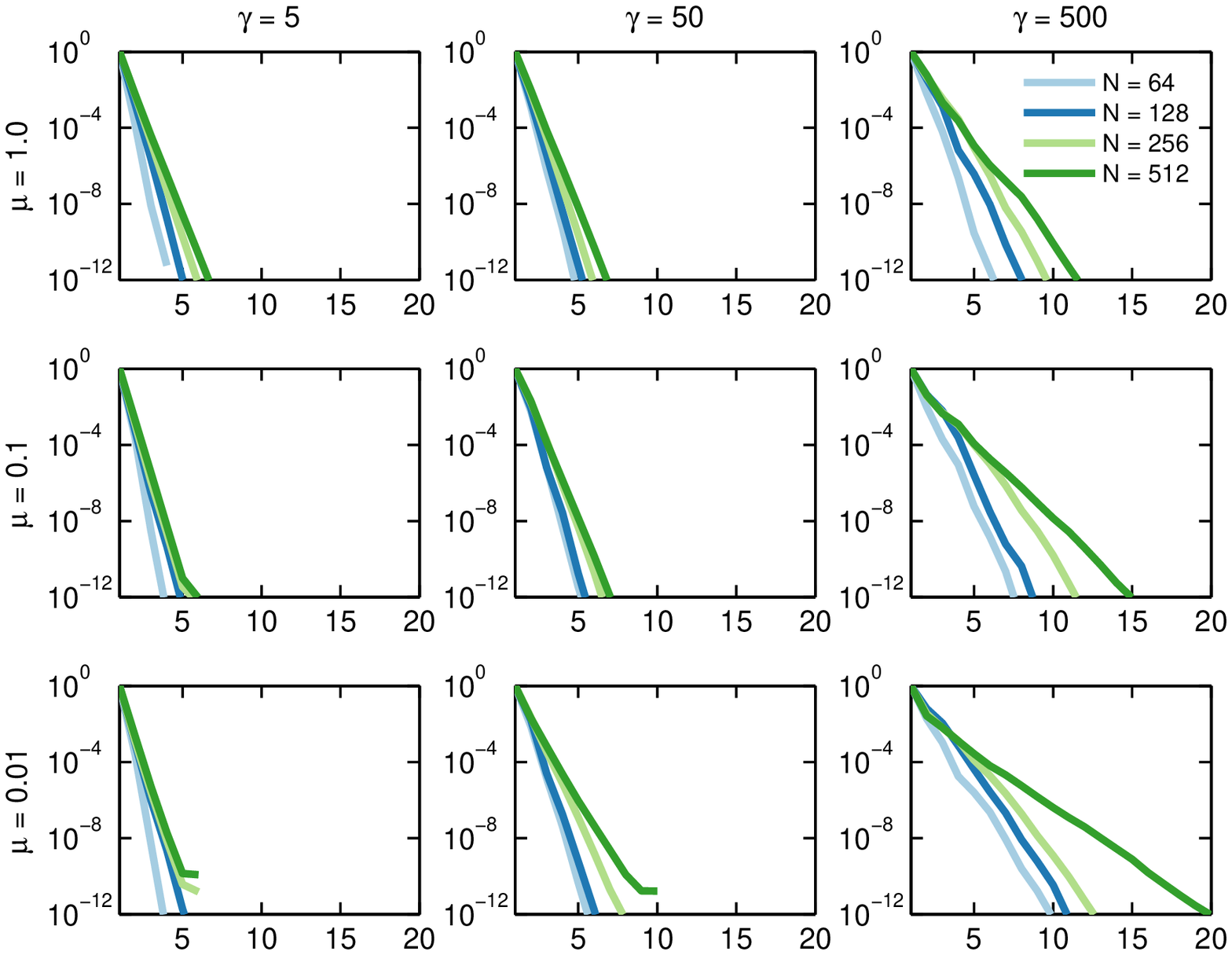}} \vspace{-1.1\baselineskip} \\
  (b) \vspace{-30pt} & \\ RAS & {\includegraphics[width=0.75\textwidth]{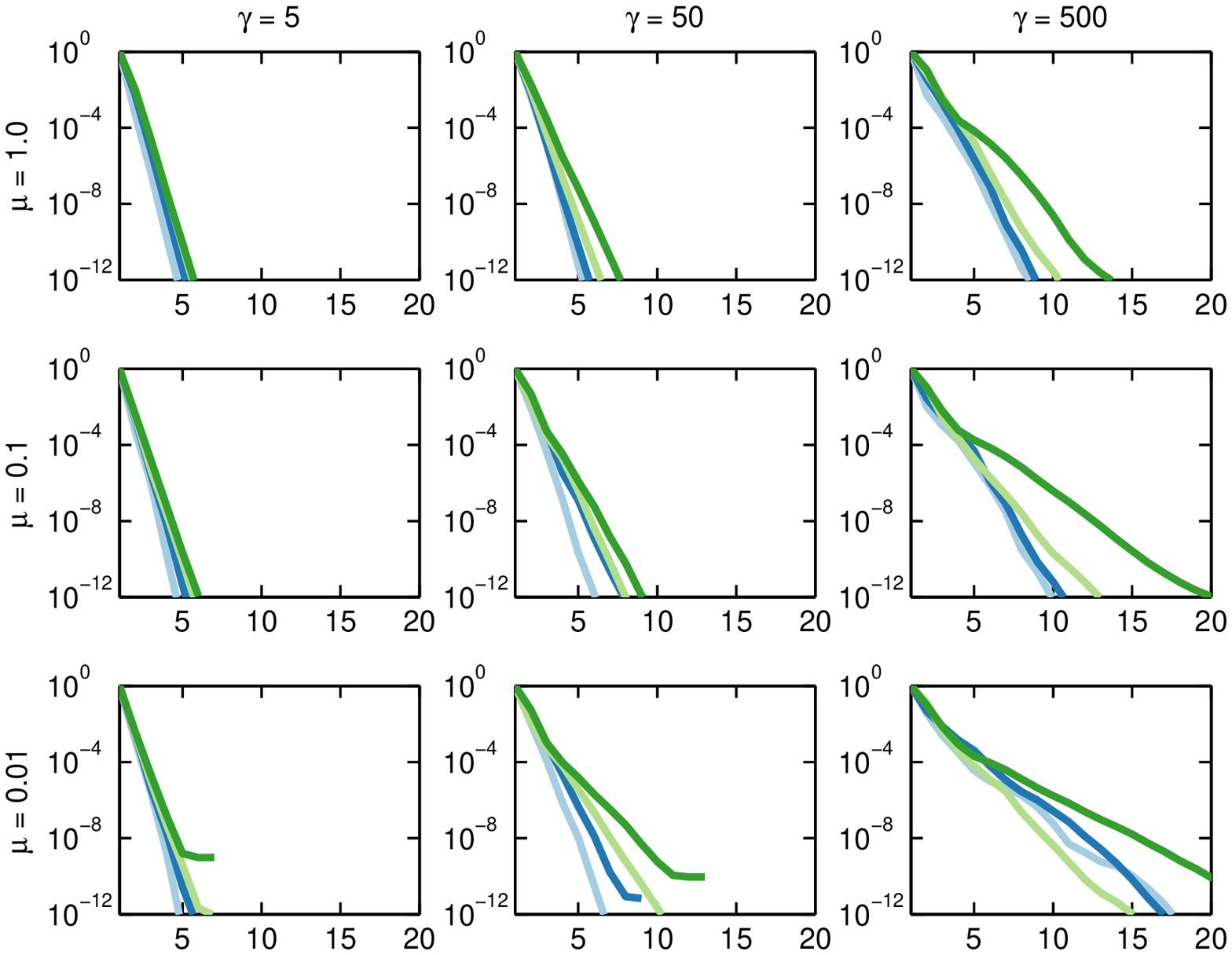}}
  \end{tabular}
\caption{
	Similar to Fig.~\ref{fig-5.1.1-3}, but here using subdomains of size $16 \times 16$ and an overlap width of $4$.
	With larger subdomains and overlap widths, the performance of the additive algorithm approaches that of the multiplicative algorithm for a broader range of physical parameters (compare to Fig.~\ref{fig-5.1.1-3}).
}
\label{fig-5.1.1-4}
\end{figure}

In these tests, we execute a single time step of the semi-implicit IB time integrator with the convective term disabled, so that we can focus on linear solver performance.
Each linear solve is allowed to run until the initial residual is reduced by $10^{-12}$ or it reaches 100 iterations.
We first examine the effect of subdomain size and overlap width on solver performance.
For these tests, we use $N = 128$ and $\gamma = 500$, which is considered ``very stiff'' in our previous work \cite{RDGuy15-gmgiib}, and consider Stokes flow (Fig.~\ref{fig-5.1.1-1}) and time-dependent flows for $\mu = 1.0$, $0.1$, $0.01$  (Fig.~\ref{fig-5.1.1-2}).
We use subdomains of size $4 \times 4$, $8 \times 8$, and $16 \times 16$ and consider $0$, $2$, or $4$ cells of overlap.
It is clear that solver performance degrades substantially as the fluid viscosity decreases.
At the highest Reynolds numbers, only the largest subdomains and overlaps yield effective solvers.
For easier cases, the multiplicative and additive solvers yield similar performance, but in many of the more challenging cases, the multiplicative smoother can converge in approximately half the iterations as the additive smoother for smaller overlap widths.
At lower Reynolds numbers and Stokes flow conditions, however, the additive and multiplicative algorithms yield similar performance.

We also examine the scalability of the multiplicative and additive smoothers under grid refinement.
Fig.~\ref{fig-5.1.1-3} shows the performance of the additive and multiplicative smoothers using $8 \times 8$ subdomains with an overlap of $2$ for various relative stiffnesses and viscosities, and Fig.~\ref{fig-5.1.1-4} shows results from similar tests using $16 \times 16$ subdomains with an overlap of $4$.
It is clear that at low Reynolds numbers or low relative stiffnesses, the multiplicative smoother yields an essentially scalable algorithm, as shown previously \cite{RDGuy15-gmgiib}.
At high elastic stiffnesses and larger Reynolds numbers, both solvers begin to break down, but with sufficiently large subdomain sizes and overlap widths, the additive version of the algorithm yields performance that is similar to that obtained by the multiplicative algorithm.

\subsubsection{Schur complement smoother}

\begin{figure}
\centering
\footnotesize
  \begin{tabular}{>{\centering\arraybackslash}m{32pt} >{\centering\arraybackslash}m{0.8\textwidth} }
  (a) \vspace{-30pt} & \\ & {\includegraphics[width=0.75\textwidth]{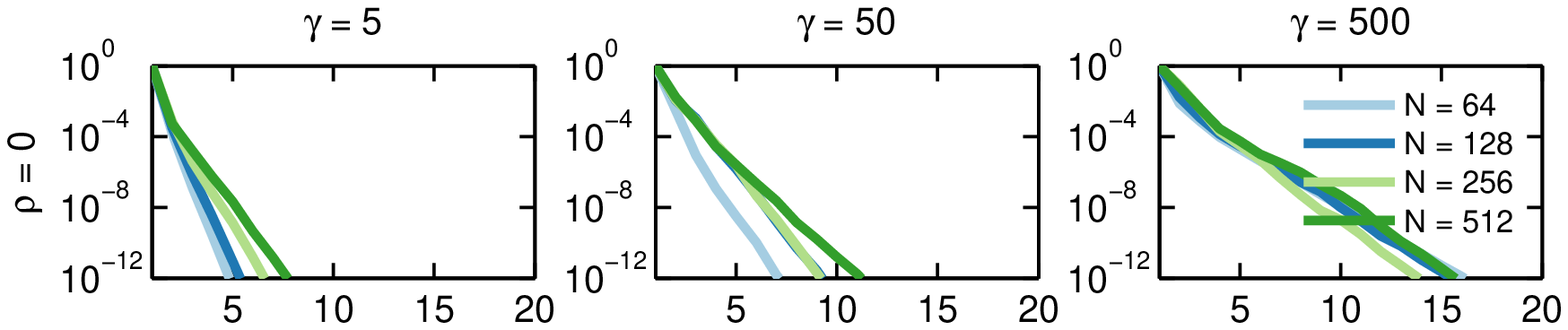}} \vspace{-1.1\baselineskip} \\
  (b) \vspace{-30pt} & \\ & {\includegraphics[width=0.75\textwidth]{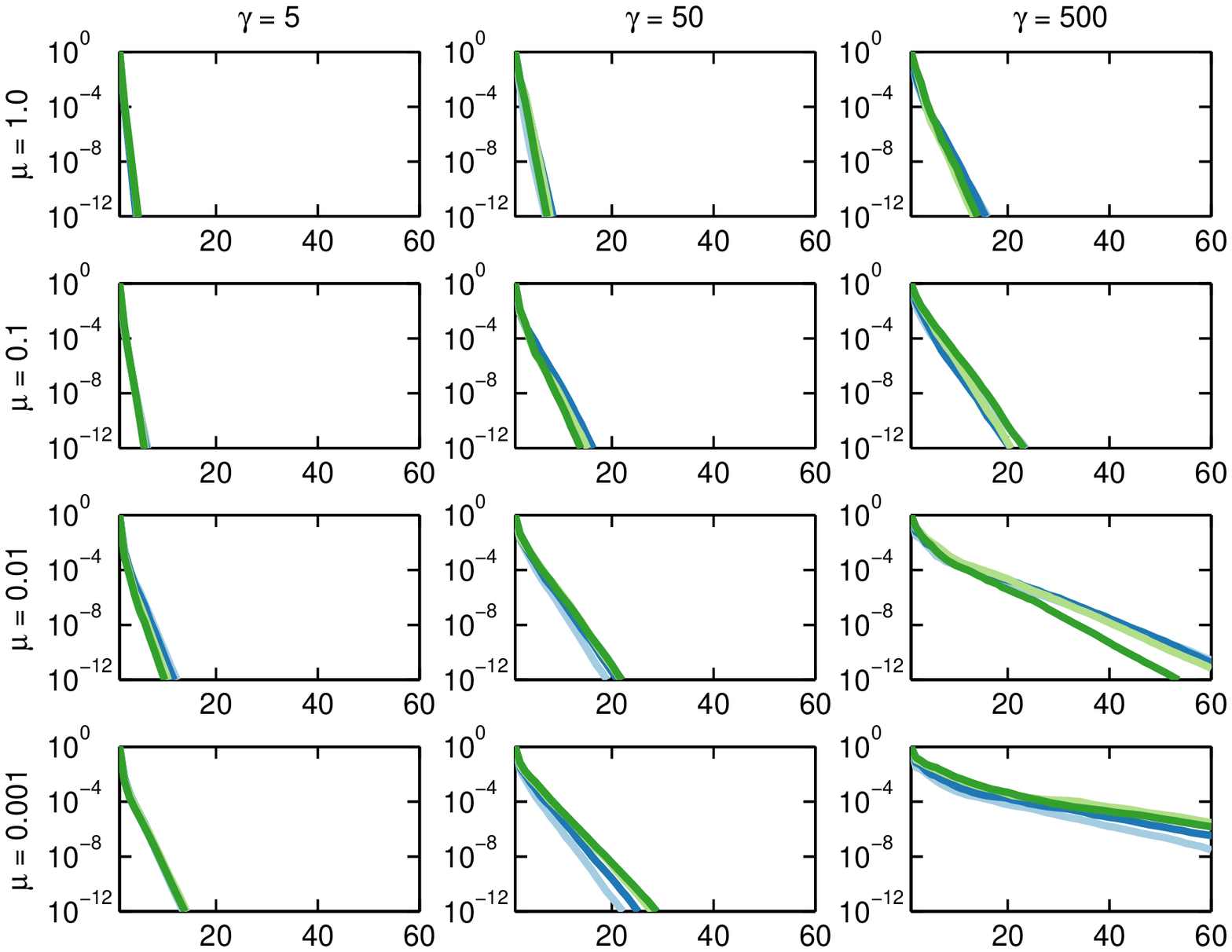}}
  \end{tabular}
\caption{
	Performance of the Schur complement-based smoother for the thick elastic shell (Sec.~\ref{sec-linear-shell}).
	We consider the effects of grid refinement for a range of relative stiffnesses ($\gamma$) and viscosities ($\mu$), for both (a)~Stokes flow conditions ($\rho = 0$) and (b)~time-dependent flow conditions with decreasing amounts of fluid viscosity.
	Notice that the Schur complement-based smoother yields a more robust algorithm except for the highest stiffnesses at low viscosities; compare to Figs~\ref{fig-5.1.1-3} and \ref{fig-5.1.1-4}.
}
\label{fig-5.1.2-1}
\end{figure}

As in the tests for the Schwarz preconditioners, we execute a single time step of the semi-implicit IB time integrator with the convective term disabled.
We perform scalability tests using the Schur complement-based smoother for Stokes flows and for time-dependent flows with various viscosities at various relative elastic stiffnesses.
Results are summarized in Fig.~\ref{fig-5.1.2-1}.
The Schur complement-based smoother is more robust under both increasing elastic stiffness and decreasing viscosity than the Schwarz-based algorithms except for the largest elastic stiffnesses.
Notice that the only case where the solver fails to reach its tight convergence threshold of $10^{-12}$ is for $\mu = 0.001$ and $\gamma = 500$.
On the other hand, the Schur complement approach generally requires somewhat more multigrid iterations than the Schwarz-based method for a given set of model parameters.
However, each application of the SC smoother is substantially less expensive than the RAS/RMS smoothers, and in our current implementation, we typically find that the SC-based solver outperforms the RAS/RMS solver in terms of total wall-clock time.

\subsection{Thin elastic membrane}
\label{sec-linear-curve}

\begin{figure}[t]
\centering
\includegraphics[width=0.625\textwidth]{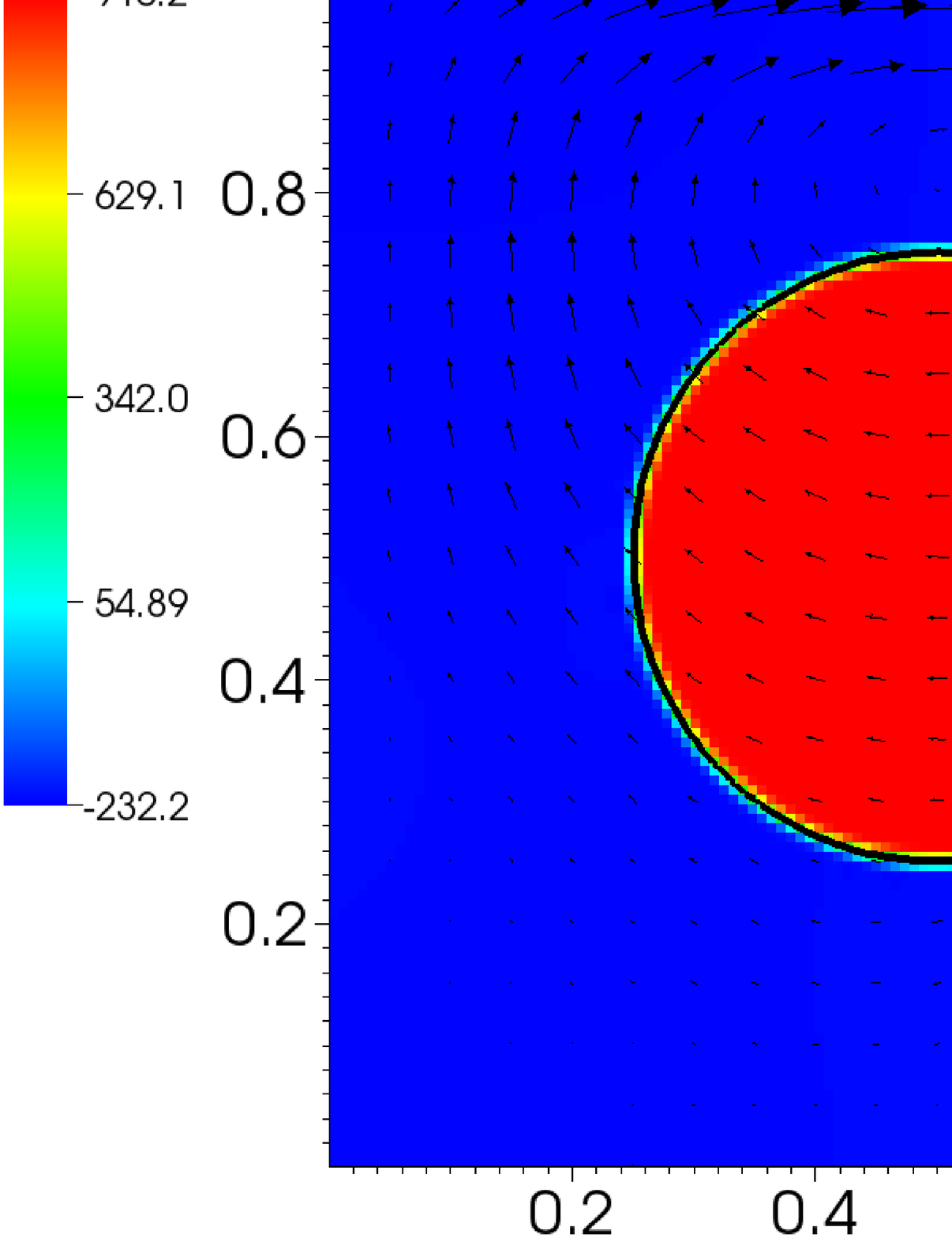}
\caption{
	Representative results from the thick elastic test problem of Sec.~\ref{sec-linear-curve}.
	Here, we plot the velocity and pressure field along with the structure configuration for $\mu=0$ and a relative stiffness of $\gamma=5$.
	The Eulerian grid spacing is $h = \frac{1}{128}$.
	Results for different values of $\mu$ and $\gamma$ are similar.
}
\label{fig-steady_curve}
\end{figure}

Next, we consider a thin circular membrane described using Lagrangian curvilinear coordinates $s_1 \in U = [0, 2 \pi)$ with initial configuration
\begin{equation}
 	\X(s_1) = \left(\xc + r\cos(s_1), \yc + r\sin(s_1) \right). \label{eqn-circle}
\end{equation}
As in Sec.~\ref{sec-linear-shell}, we choose the center to be $\x_\text{c} = (\xc, \yc) = (0.5, 0.5)$ and the radius to be $r = 1/4$.
We again use
\begin{equation}
 \F = \alpha \DD{\X}{s_1}, \label{eqn-F-lin-functional}
\end{equation}
and, as before, the Eulerian domain is discretized using an $N \times N$ grid, and the Lagrangian coordinates are discretized using $M_1 = \frac{19}{8} N$ points in the $s_1$ direction.

We again characterize the stiffness $\alpha$ in terms of a stiffness ratio $\gamma$, but $\alpha$ is now defined via
\begin{equation}
	\alpha = 7 \gamma \frac{3.93}{0.005},
\end{equation}
which yields approximately the same total force as in the thick interface case.
In this thin case, at a fixed time step size, the maximum stiffness allowed by an explicit solver \emph{decreases} in proportion to the grid spacing.
This is in contrast to the thick case.
Thus, for a fixed mechanical stiffness, the \emph{numerical} stiffness of the problem increases under grid refinement.
Consequently, the thin interface case poses substantially greater challenges to the solvers.
Moreover, an analysis similar to that presented for a thick elastic shell \cite{RDGuy15-gmgiib} implies that $\gamma = 5$ is approximately a factor of 100 times stiffer than the largest elastic stiffness permitted by an explicit solver at $\dx = \frac{1}{64}$ and $\dt = 0.005$.
Thus, the thin cases considered here are much more numerically challenging than the thick cases considered above.
Representative results are shown in Fig.~\ref{fig-steady_curve}.

\subsubsection{Schwarz smoothers}

\begin{figure}[t]
\centering
\includegraphics[width=0.6\textwidth]{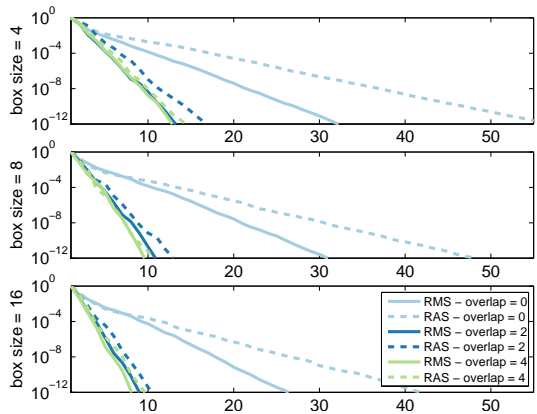}
\caption{
	Multiplicative (solid line) and additive (dashed line) smoother performance under Stokes flow conditions for the thin elastic membrane (Sec.~\ref{sec-linear-curve}) with $N = 128$ and $\gamma = 500$, and for different subdomain sizes overlap widths.
	Unlike the case of a thick elastic shell, in this case solver performance has a strong dependence on overlap width.
	With an overlap of 4, the additive and multiplicative smoothers yield similar performance.
}
\label{fig-5.2.1-1}
\end{figure}

\begin{figure}
\centering
\footnotesize
  \begin{tabular}{>{\centering\arraybackslash}m{32pt} >{\centering\arraybackslash}m{0.8\textwidth} }
  (a) \vspace{-30pt} & \\ $\mu = 1.0$   & {\includegraphics[width=0.6\textwidth]{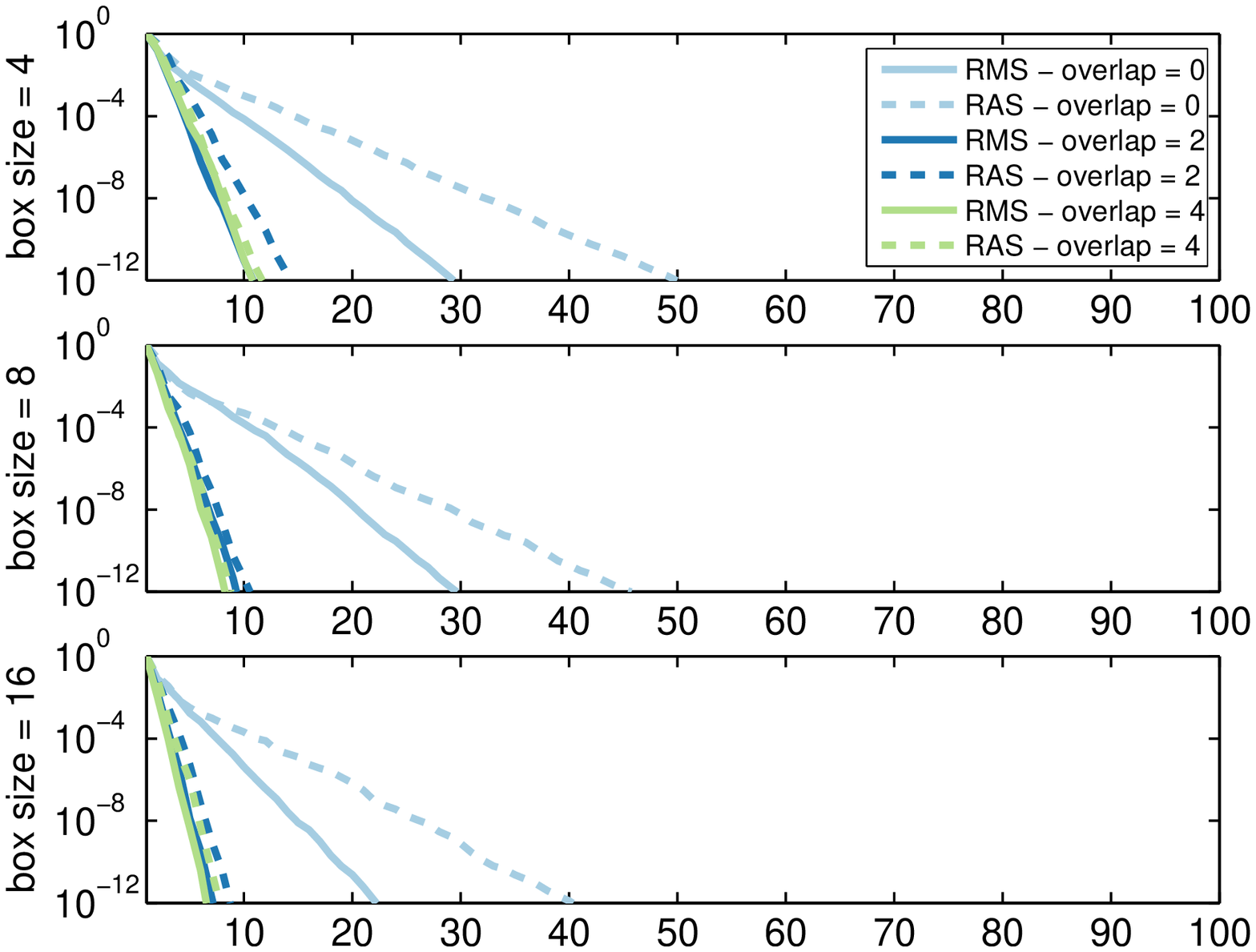}} \vspace{-1.1\baselineskip} \\
  (b) \vspace{-30pt} & \\ $\mu = 0.1$   & {\includegraphics[width=0.6\textwidth]{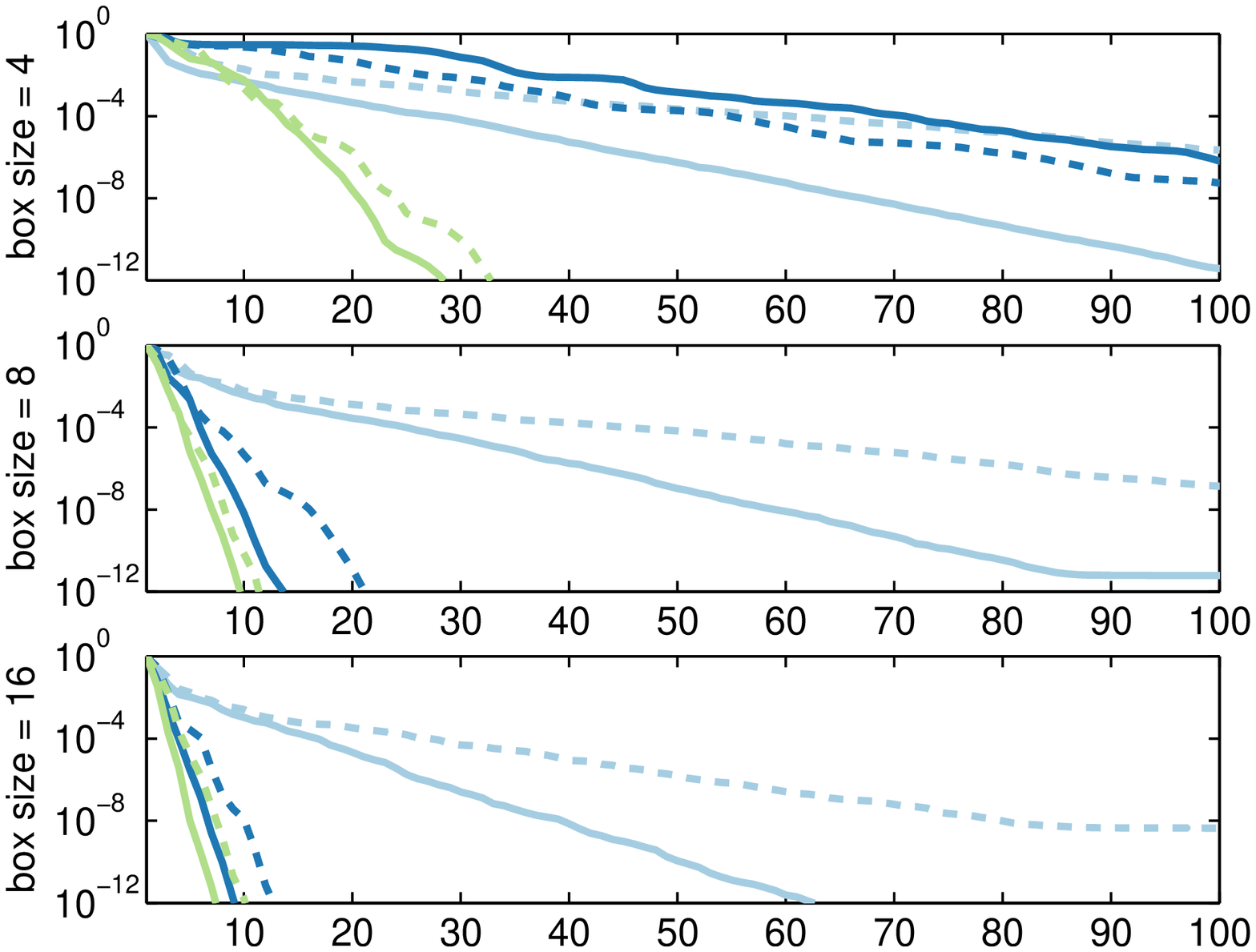}} \vspace{-1.1\baselineskip} \\
  (c) \vspace{-30pt} & \\ $\mu = 0.01$  & {\includegraphics[width=0.6\textwidth]{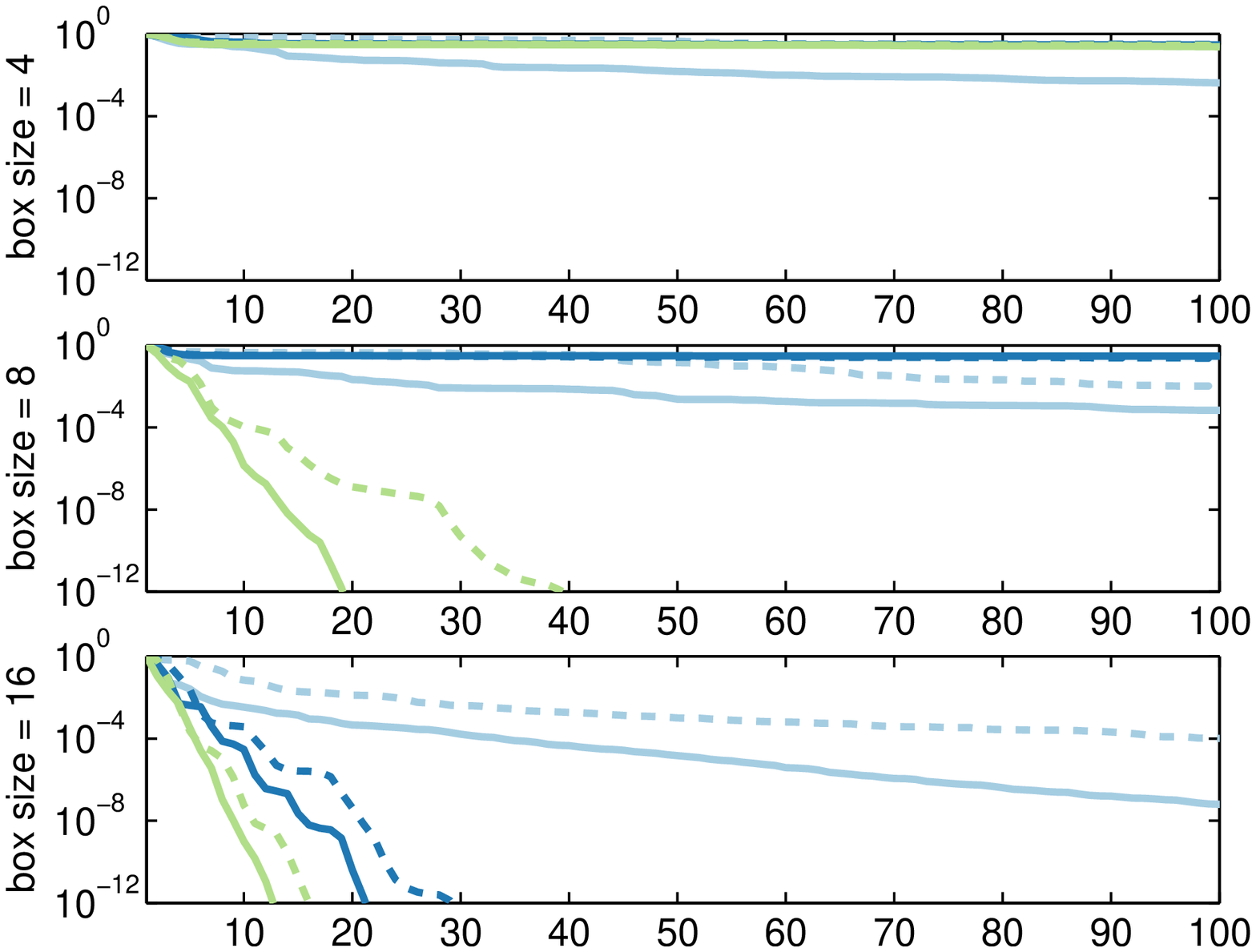}}
  \end{tabular}
\caption{
	Similar to Fig.~\ref{fig-5.2.1-1}, but here showing the effect of decreasing viscosity, with (a)~$\mu = 1$, (b)~$\mu = 0.1$, and (c)~$\mu = 0.01$.
	As in the thick case (Sec.~\ref{sec-linear-shell}), performance degrades with decreasing viscosity, and the multiplicative smoother can outperform the additive smoother by a wide margin.
}
\label{fig-5.2.1-2}
\end{figure}

\begin{figure}
\centering
\footnotesize
  \begin{tabular}{>{\centering\arraybackslash}m{32pt} >{\centering\arraybackslash}m{0.8\textwidth} }
  (a) \vspace{-30pt} & \\ RMS & {\includegraphics[width=0.75\textwidth]{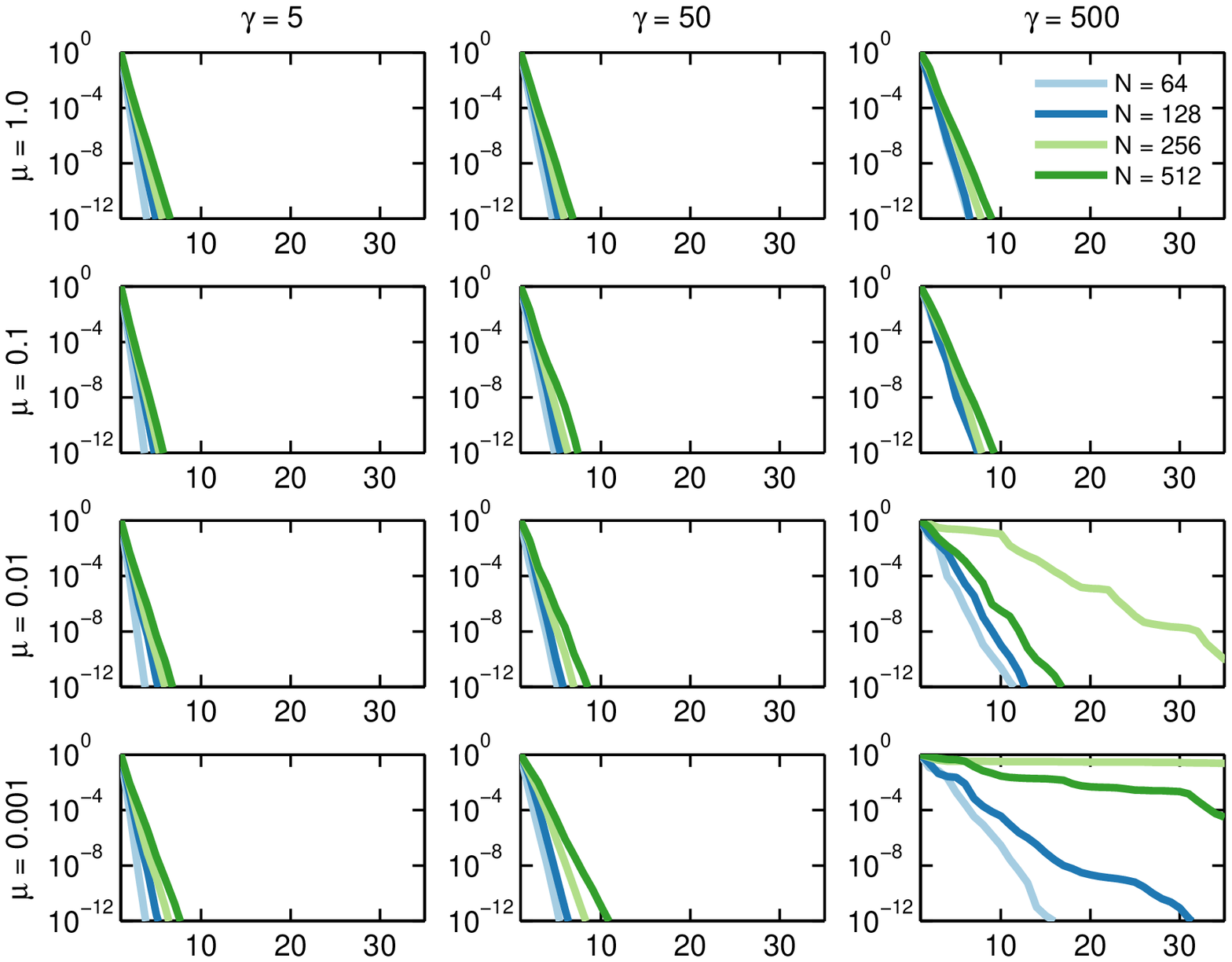}} \vspace{-1.1\baselineskip} \\
  (b) \vspace{-30pt} & \\ RAS & {\includegraphics[width=0.75\textwidth]{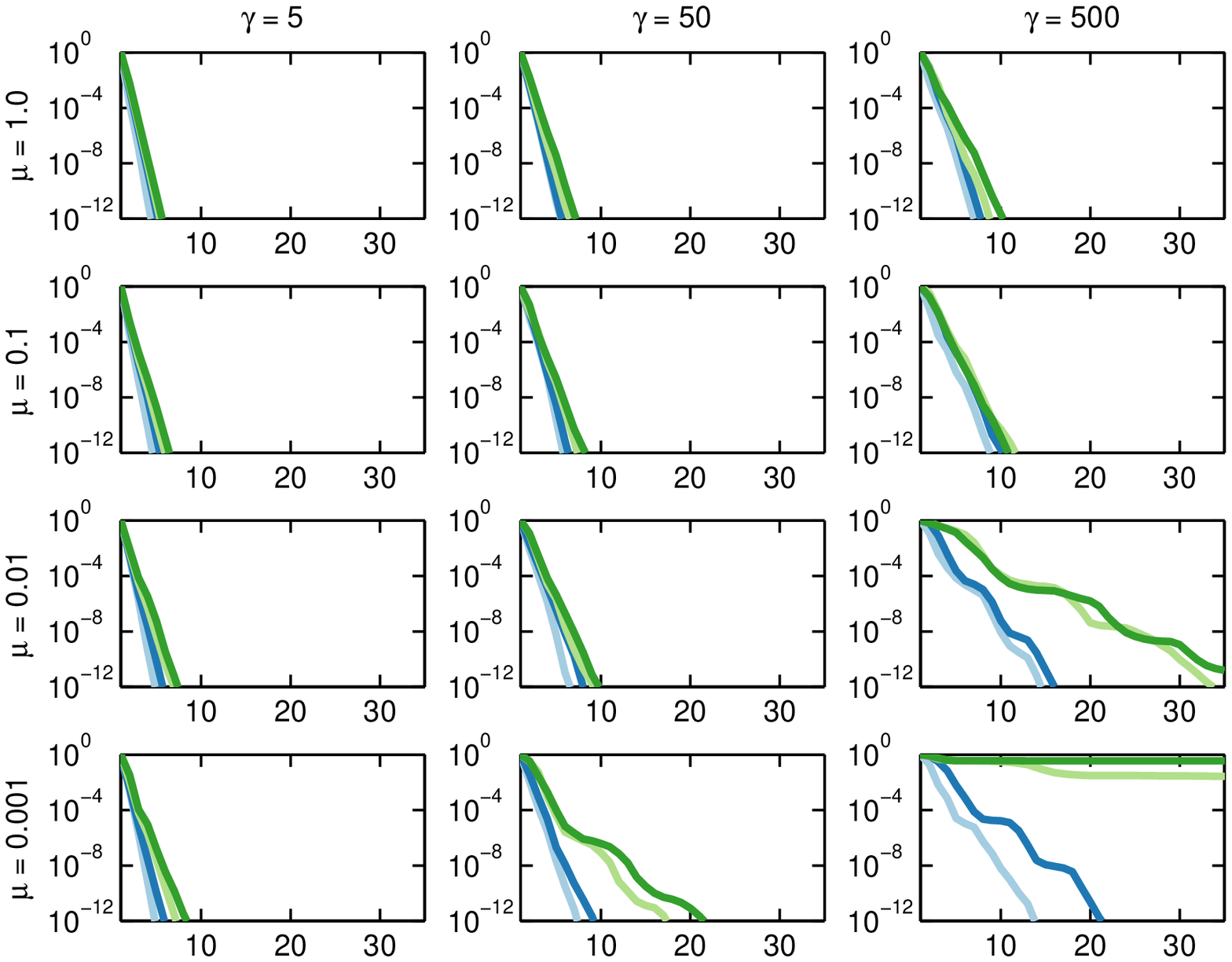}}
  \end{tabular}
\caption{
	Performance of multigrid using (a) multiplicative (RMS) and (b) additive (RAS) smoothers under grid refinement for a range of relative stiffnesses ($\gamma$) and viscosities ($\mu$), using subdomains of size $16 \times 16$ and an overlap width of $4$.
	Notice that the performance of the additive algorithm is similar to that of the multiplicative algorithm in most cases.
	Both algorithms ultimately stall for the largest stiffnesses for sufficiently small viscosity.
}
\label{fig-5.2.1-3}
\end{figure}

As before, we disable the convective term in our semi-implicit time integrator to focus on linear solver performance.
Fig.~\ref{fig-5.2.1-1} shows the effect of subdomain size and overlap width on solver performance for the multiplicative and additive smoothers in Stokes flow conditions at a relative stiffness of $\gamma = 500$.
In this thin interface case, solver performance depends strongly on overlap width.
This is in contrast to the case of a thick elastic shell (e.g.~Fig.~\ref{fig-5.1.1-1}).
At an overlap width of $4$, the two Schwarz smoothers yield nearly identical convergence rates for all subdomain sizes considered in the Stokes flow case.
Time dependent flow conditions pose a greater challenge to the solver (Fig.~\ref{fig-5.2.1-2}).
At $\mu = 1$, solver performance is largely insensitive to subdomain size and overlap width, except for an overlap width of $0$.
By contrast, for $\mu = 0.1$, there are substantial differences in solver performance for the different subdomain sizes and overlap widths.
For $\mu = 0.01$, only the largest subdomain sizes and overlap widths yield convergent solver algorithms.
It is clear that the thin interface case is fundamentally \emph{more stiff} than the thick body case.

Fig.~\ref{fig-5.2.1-3} shows the effects of grid refinement on solver performance for the time-dependent cases.
For these tests, we consider only a subdomain size of $16 \times 16$ along with an overlap width of $4$.
The multiplicative smoother yields an essentially scalable multigrid algorithm except for the highest stiffness ($\gamma = 500$) and lowest viscosity ($\mu = 0.01$) considered.
Performance of the additive smoother is similar except for $\gamma = 500$ and $\mu = 0.01$.
In this challenging case, both smoothers show poor performance for $N = 256$, and the additive smoother stagnates at the highest grid spacing ($N = 512$).
In practical time-dependent calculations, we likely would use a relative convergence tolerance around $10^{-6}$, which corresponds to 4--5 multigrid iterations in all but the most difficult cases considered here.

\subsubsection{Schur complement smoothers}

\begin{figure}
\centering
\footnotesize
  \begin{tabular}{>{\centering\arraybackslash}m{32pt} >{\centering\arraybackslash}m{0.8\textwidth} }
  (a) \vspace{-30pt} & \\ & {\includegraphics[width=0.75\textwidth]{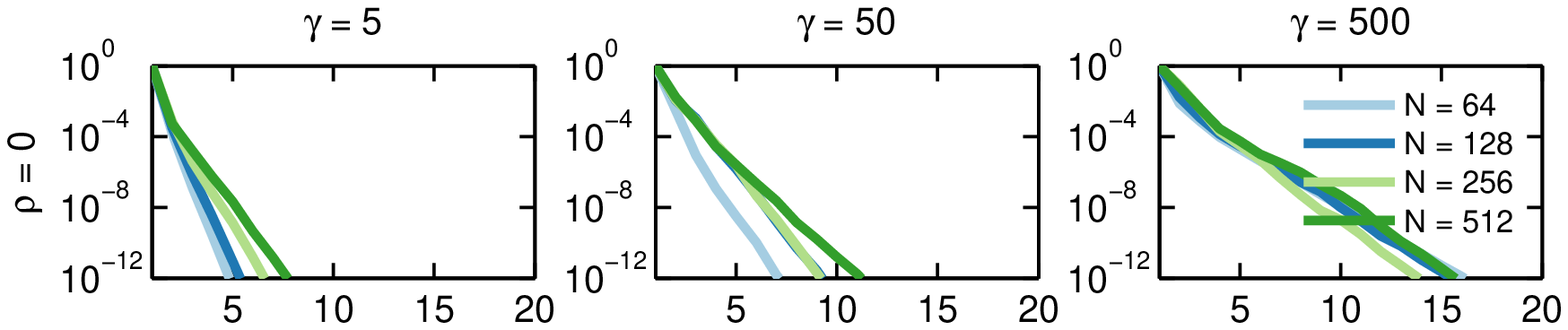}} \vspace{-1.1\baselineskip} \\
  (b) \vspace{-30pt} & \\ & {\includegraphics[width=0.75\textwidth]{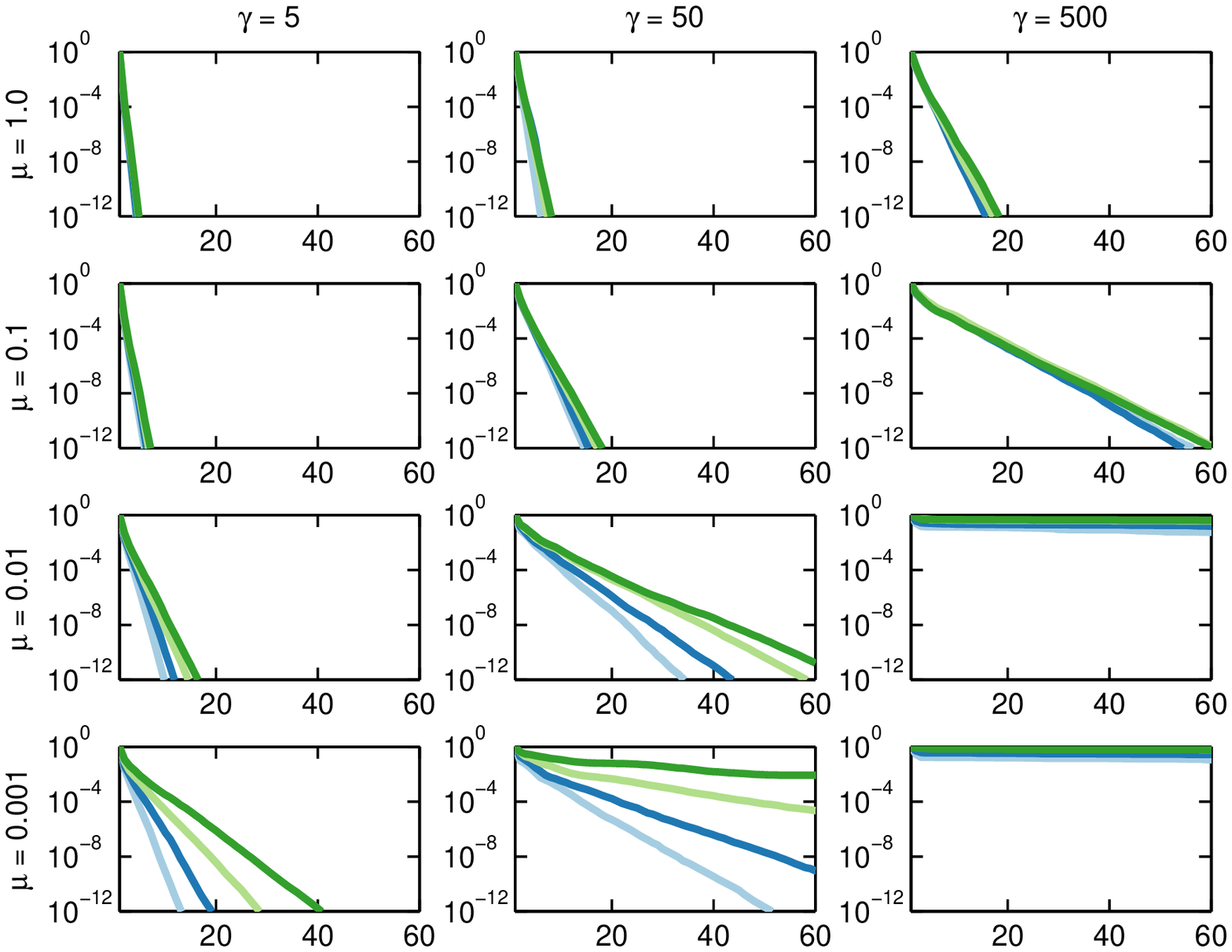}}
  \end{tabular}
\caption{
	Performance of the Schur complement-based smoother for the thin elastic membrane (Sec.~\ref{sec-linear-curve}).
	We consider the effects of grid refinement for a range of relative stiffnesses ($\gamma$) and viscosities ($\mu$), for (a)~Stokes flow conditions ($\rho = 0$) and (b)~time-dependent flow conditions with decreasing amounts of fluid viscosity.
	The Schur complement-based smoother is very effective for Stokes flow conditions and lower viscosities, but at higher Reynolds numbers, performance degrades even at relatively small stiffnesses ($\gamma = 5$), and the solver stagnates for the higher stiffnesses considered.
}
\label{fig-5.2.2-1}
\end{figure}

As in the tests for the Schwarz preconditioners, we execute a single time step of the semi-implicit IB time integrator with the convective term disabled, now using the Schur complement-based smoother.
Fig.~\ref{fig-5.2.2-1} summarizes solver performance under grid refinement for a range of flow conditions and elastic stiffnesses.
The Schur complement-based smoother is extremely robust for Stokes flows and low Reynolds number cases.
However, the solver begins to stagnate for the higher stiffness cases even at a modest Reynolds number of 100.
For this thin interface case, it appears that the current additive and multiplicative Schwarz smoothers are more effective, although they also struggle with higher Reynolds numbers and elastic stiffnesses.

\subsection{Suspension of immersed structures}
\label{sec-linear-suspension}

\begin{figure}[t]
\centering
\footnotesize
  \begin{tabular}{cccc}
  (a) & & (b) & \\
  & \includegraphics[width=0.4\textwidth]{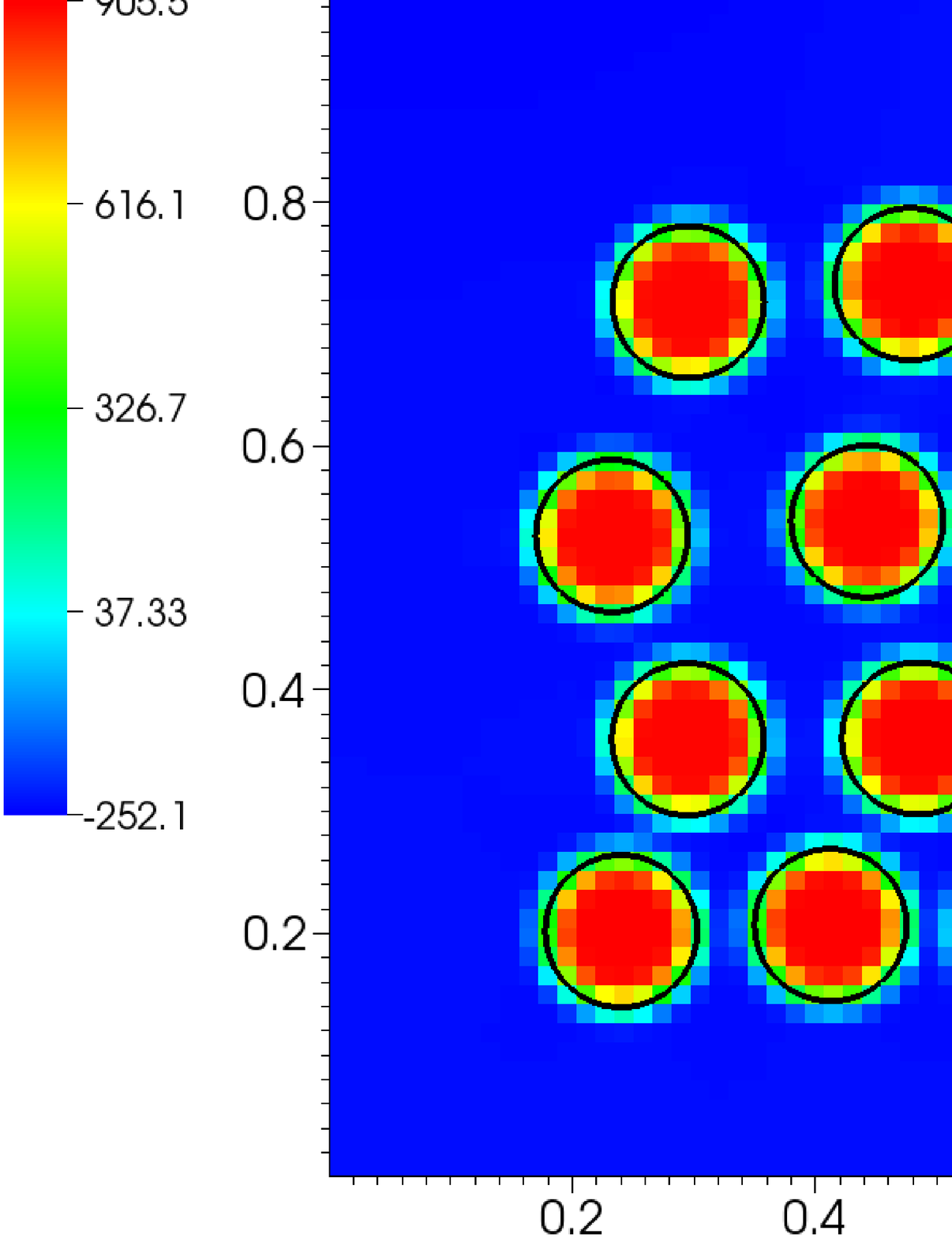} &
  & \includegraphics[width=0.4\textwidth]{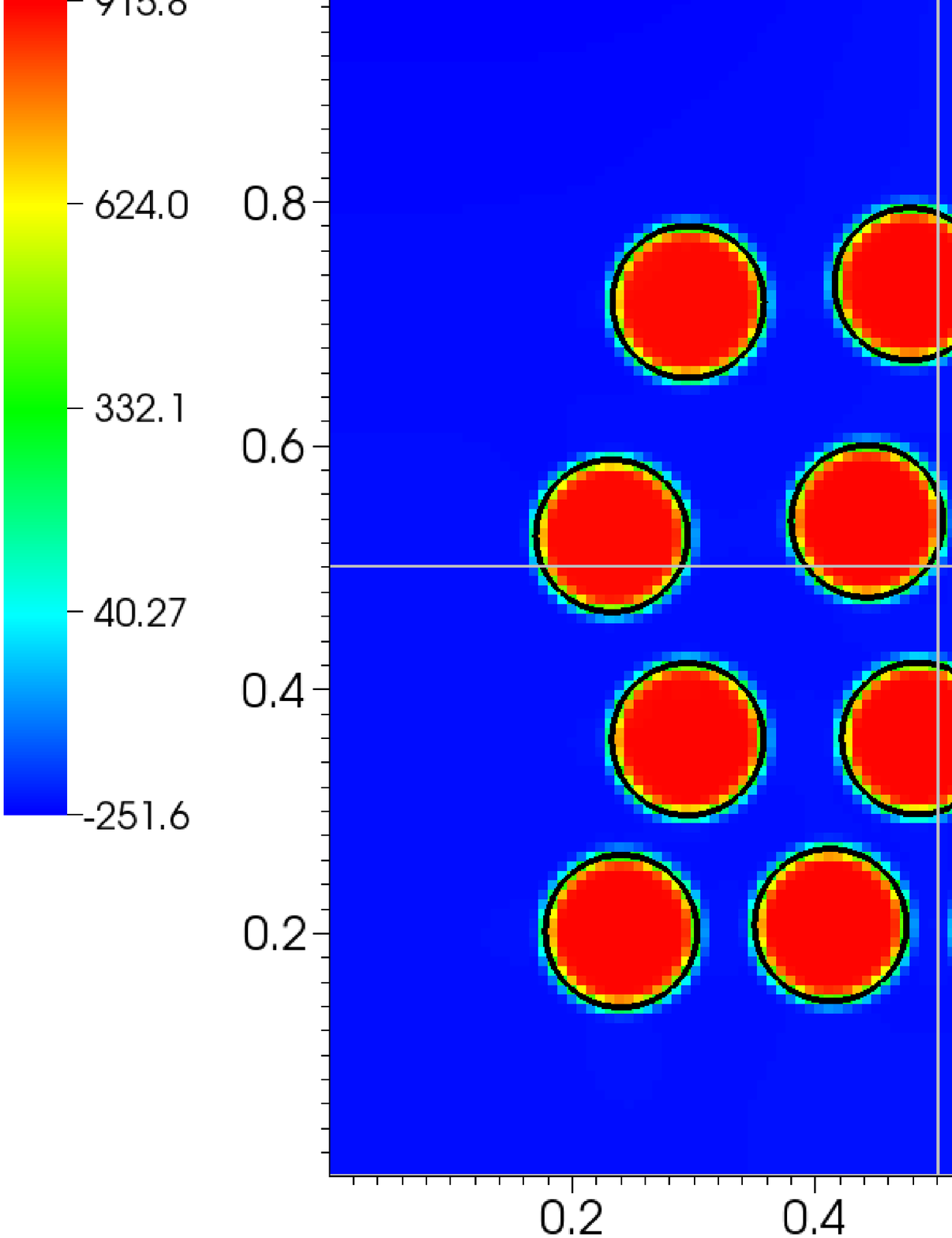} \\
  (c) & & (d) & \\
  & \includegraphics[width=0.4\textwidth]{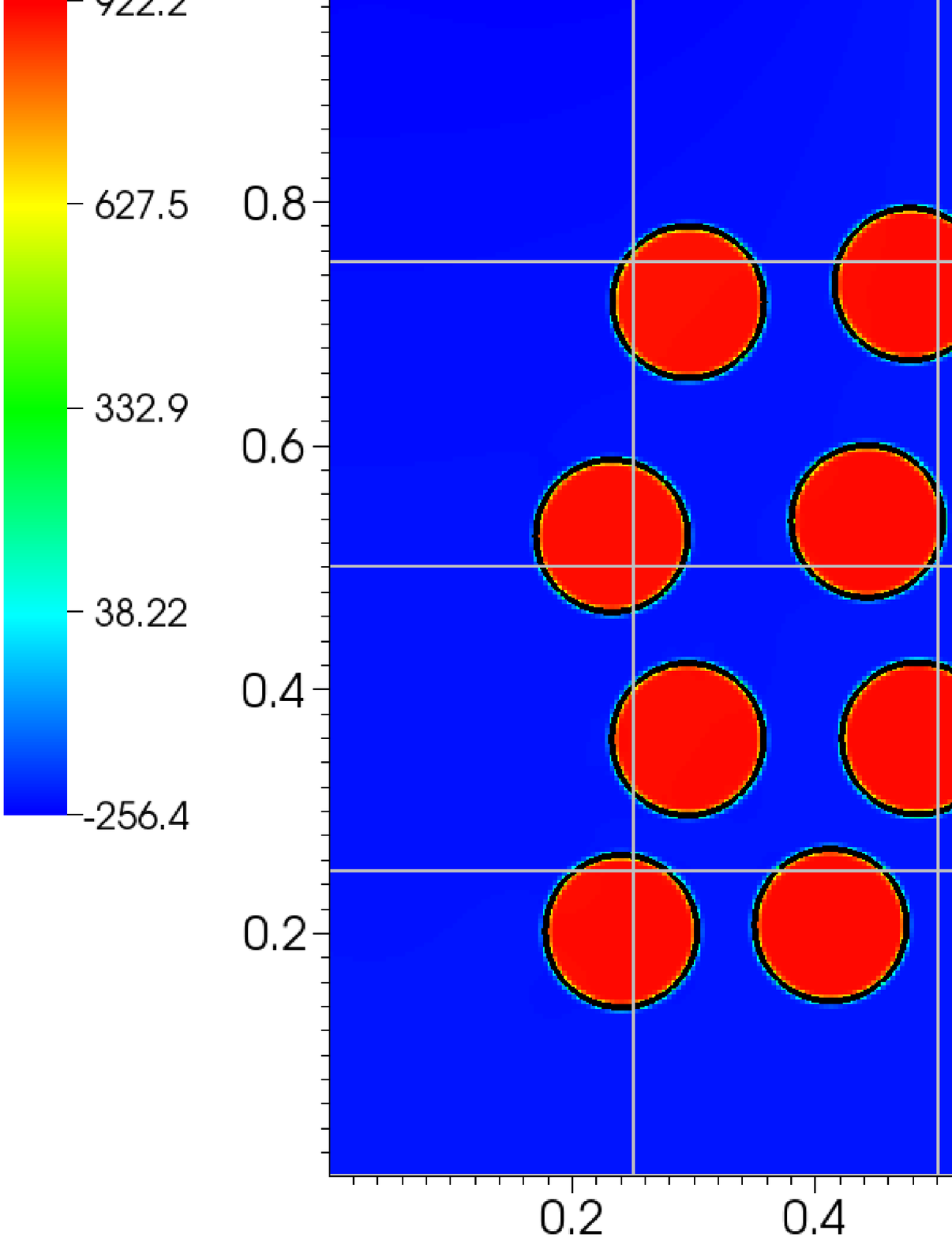} &
  & \includegraphics[width=0.4\textwidth]{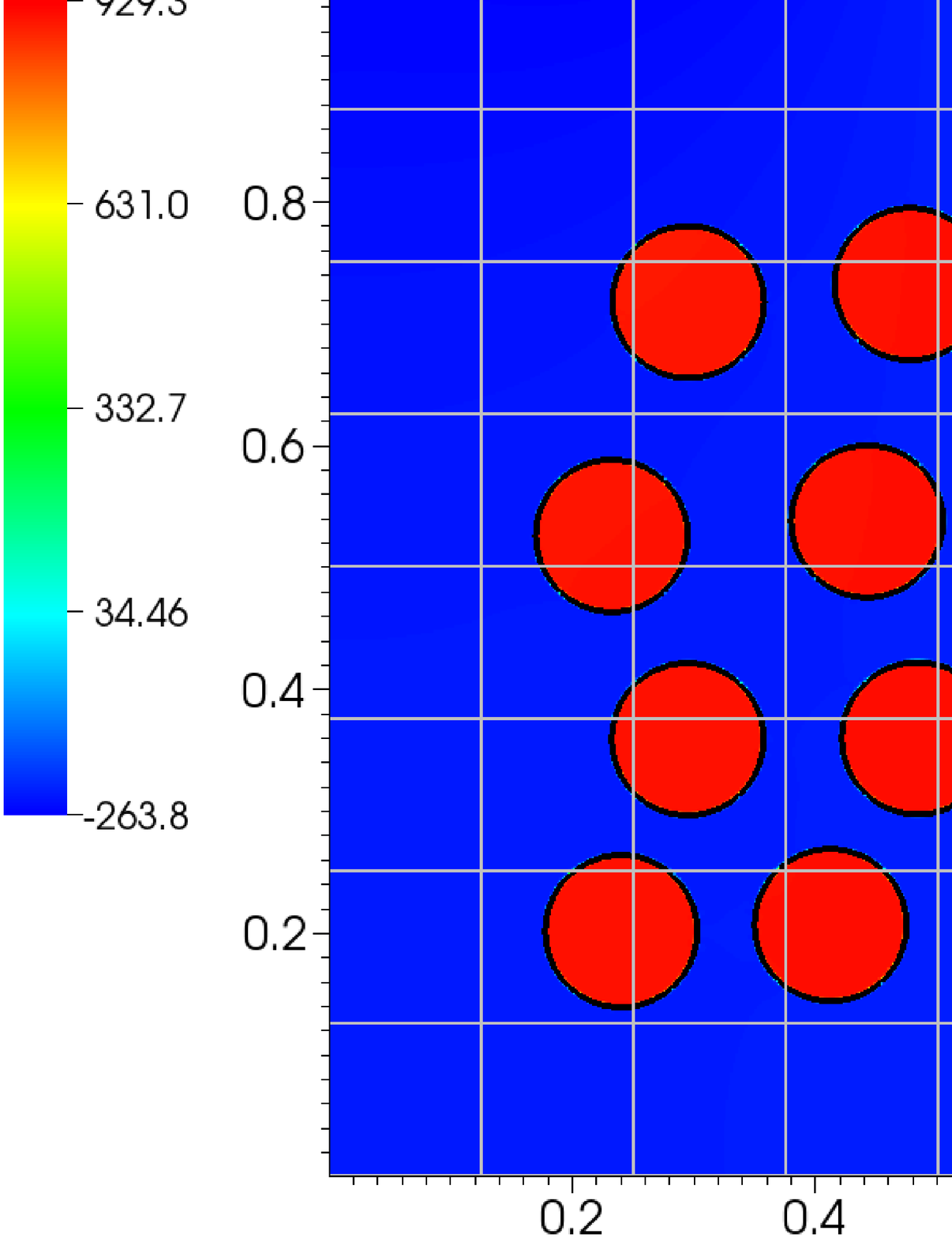} \\
  \end{tabular}
\caption{
  Pressure field and parallel domain decomposition (indicated by gray boxes) for a suspension of 16 circular immersed interfaces with $\rho = 1$, $\mu = 1$, and $\gamma = 5$ for (a)~$N=64$, (b)~$N=128$, (c)~$N=256$, and (d)~$N=512$.
}
\label{fig-5.3-1}
\end{figure}

\begin{figure}
\centering
\footnotesize
  \begin{tabular}{>{\centering\arraybackslash}m{32pt} >{\centering\arraybackslash}m{0.8\textwidth} }
  (a) \vspace{-30pt} & \\ RAS & {\includegraphics[width=0.8\textwidth]{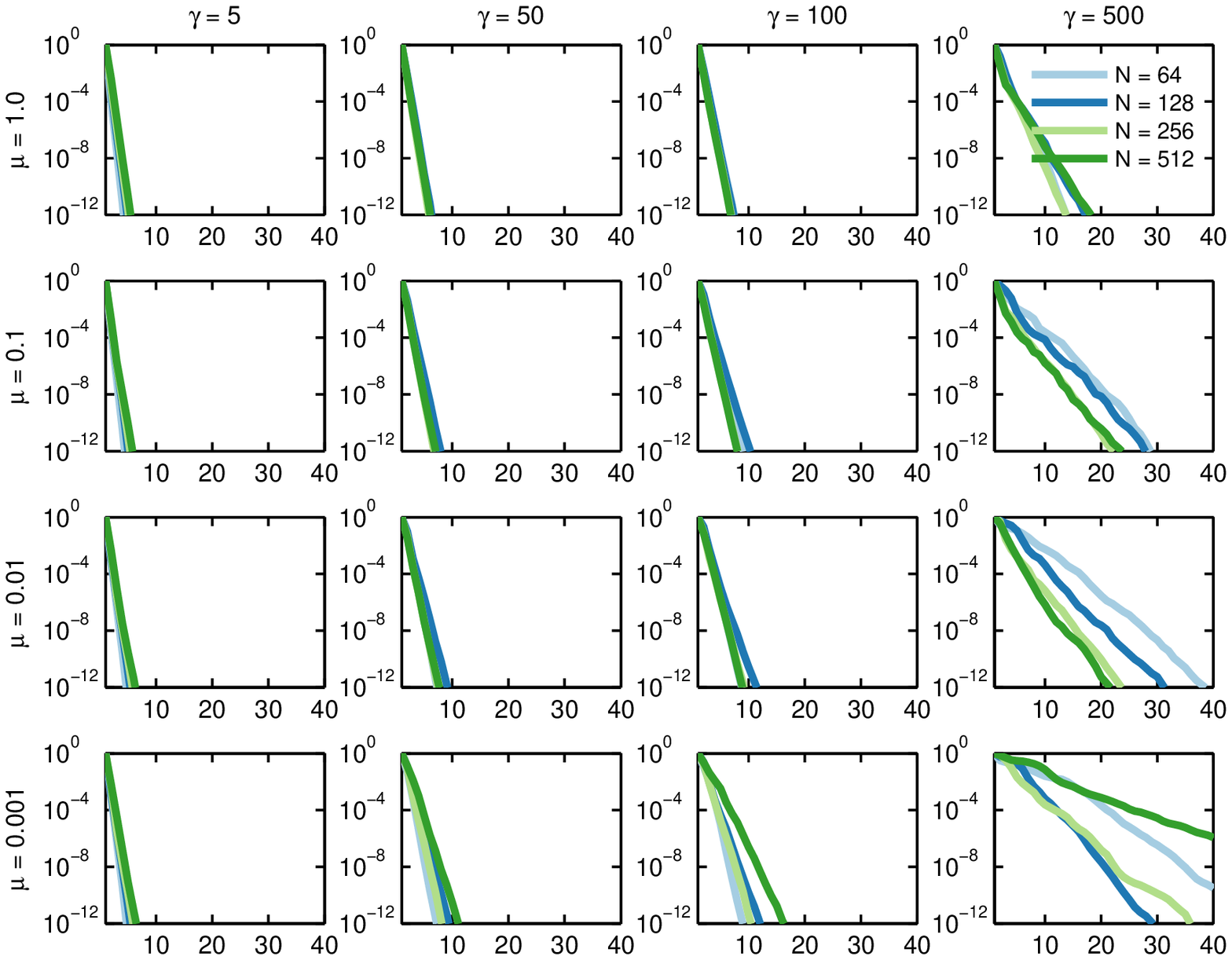}} \vspace{-2\baselineskip} \\
  (b) \vspace{-30pt} & \\ SC  & {\includegraphics[width=0.8\textwidth]{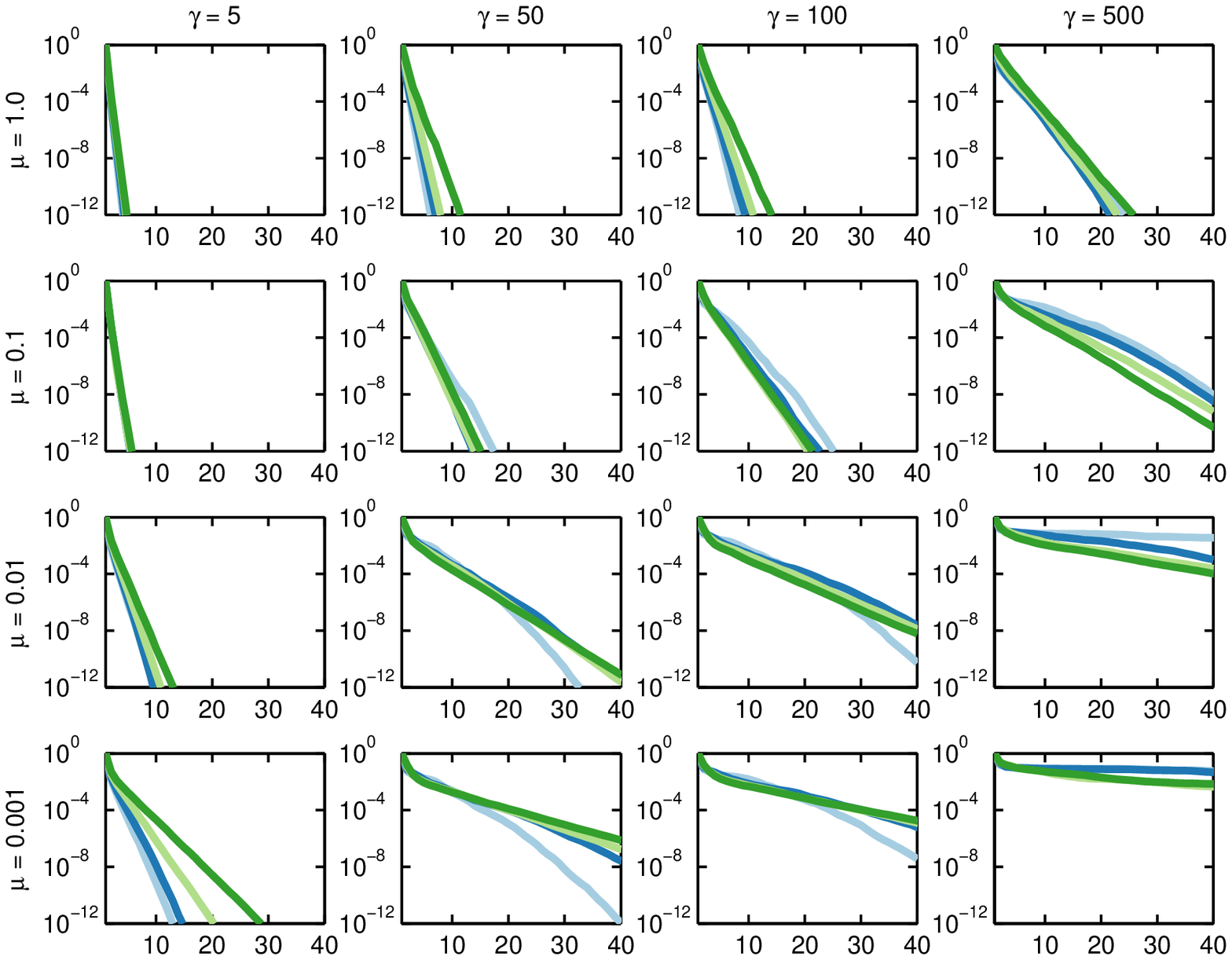}}
  \end{tabular}
\caption{
	Parallel scalability using the (a)~RAS smoother and (b)~Schur complement smoother at nonzero Reynolds number flow conditions for a suspension of elastic membranes (Sec.~\ref{sec-linear-suspension}) using subdomains of size $16 \times 16$ with overlap widths of $4$.
	We use $(N/64)^2$ processors for each case, so that the number of grid cells assigned to each processor remains fixed.
	Both solvers are essentially scalable for the lower Reynolds number cases, but as in the serial case, performance degrades with increasing elastic stiffness ($\gamma$) and decreasing viscosity ($\mu$).
}
\label{fig-5.3-2}
\end{figure}

This test case is similar to the thin interface case of Sec.~\ref{sec-linear-curve}, but here we consider a suspension of 16 structures, each with an initial configuration corresponding to a circle of radius $r = 1/16$.
The structures are randomly placed in the domain and are required not to overlap each other or the domain boundary.
Fig.~\ref{fig-5.3-1} shows the distribution of structures along with the resulting pressure field for $\gamma = 5$.
We consider only the RAS and Schur complement smoothers, and we explore the performance of the solver with increasing numbers of processors for $N = 64$, $128$, $256$, and $512$, using $(N/64)^2$ processors for each case, so that the number of grid cells assigned to each processor remains fixed.
Fig.~\ref{fig-5.3-1} shows the Cartesian grid-based parallel domain decompositions.
Fig.~\ref{fig-5.3-2} summarizes the solver performance under grid refinement for a range of flow conditions and elastic stiffnesses.
Performance is similar to that obtained in serial for the case of a single immersed membrane, although the Schur complement-based smoother shows slightly poorer scaling in parallel than in serial.
This is not unexpected because the Schur complement-based solver uses processor-restricted Gauss-Seidel in its subdomain operators rather than a true parallel Gauss-Seidel algorithm.
Nonetheless, the Schur complement-based solver yields good scalability in cases where the underlying serial algorithm also yields good scalability.
As also observed in the serial case, the Schur complement-based algorithm ultimately stalls for sufficiently small viscosities or sufficiently large elastic stiffnesses.

%

\section{Discussion and conclusions}

This paper has extended a geometric multigrid (GMG) preconditioning approach to semi-implicit formulations of the immersed boundary (IB) method~\cite{RDGuy12,RDGuy15-gmgiib} in several important ways.
First, we showed that the multiplicative ``big-box'' Vanka smoother previously developed by Guy et al.~\cite{RDGuy15-gmgiib} can be recast as Richardson iterations preconditioned by a multiplicative Schwarz domain decomposition method, and we demonstrated that a restricted additive Schwarz (RAS)~\cite{XCCai99, EEfstathiou03-RAS} variant of this algorithm is also an effective smoother for the Stokes-IB systems that occur in this semi-implicit formulation.
Although RAS yields convergence rates that are lower than multiplicative Schwarz, extending the smoother to a purely additive algorithm is crucial for deploying these methods in parallel computing environments because multiplicative domain decomposition methods impose a sequentiality that is not amenable to large-scale parallelization.
Indeed, in the limit of large numbers of processors, it is clear that an additive smoother that does not require or assume a particular order in which the subdomains is processed is essential to achieving good parallel scalability.

We further demonstrated that an even more effective smoother approach is obtained by considering an approximate block factorization of the Stokes-IB operator that appears in our semi-implicit formulation.
What is remarkable about this Schur complement-based smoother is that it requires only a few iterations of point-relaxation smoothers on suitably constructed block operators for the velocity and pressure degrees of freedom.
Consequently, the computational complexity of a single application of this smoother is comparable in complexity to optimal multigrid smoothers for much simpler systems such as isotropic Poisson problems.
The SC smoother is also additive and well-suited for large-scale parallelization.
Similar Stokes-type operator also appear in geodynamic applications that consider strong anisotropic viscosity variations (e.g.~in the work of May et al.~\cite{May15-geodynamics_mg_fem,May08-geodynamics_mg_pc,Furuichi11-geodynamics_schur}, which has proposed scalable multigrid preconditioners for such applications), and our Schur complement is similar to those used in this earlier work.

We performed extensive tests of the GMG algorithm using both RAS/RMS and SC smoothers.
As in earlier work \cite{RDGuy15-gmgiib}, we observe that solver performance degrades with increasing elastic stiffness.
This study also reveals that the present solver approach degrades with increasing Reynolds number, with all methods ultimately failing for sufficiently small fluid viscosities.
On the other hand, both the RAS/RMS and SC smoothers were shown to yield nearly optimal convergence rates at low Reynolds numbers and in Stokes flow conditions.
Consequently, these methods may ultimately prove to offer practical solver strategies for important biological applications at the cellular and sub-cellular scales.
The extension of this methodology to moderate-to-high Reynolds numbers remains important future work.

Although the present study considers only linear solver performance, this linear solver is implemented within a time stepping framework that supports both linear and nonlinear structural models.
In the nonlinear case, we use a Newton-Krylov method \cite{KnollKeyes04}, which requires solutions to systems of the form~\eqref{eqn-imp-momentum}--\eqref{eqn-imp-F}, but with $\cF_h$ replaced by a linearized force operator $\cK_h = \left.\left(\D{\cF_h}{\X}\right)\right|_{\X = \X^{n+1}}$ for successive approximations to $\X^{n+1}$.
We have found that because the configuration of the structure does not change very much within a time step, we generally can successfully use the configuration $\X^{n}$ to construct a ``lagged'' preconditioner.
JFNK generally appears to be quite effective so long as the underlying linear solver algorithm is effective.
At present, however, the performance of both the linear and nonlinear implicit time stepping schemes lags that of our more mature explicit dynamics codes.
In Stokes flow conditions and at very high stiffness, the implicit solver can yield wall-clock times comparable to our explicit solvers.
In most other cases, however, the implicit solver generally remains a factor of 2--10 slower than the explicit solver, despite the fact that the implicit solver is able to use much larger time step sizes than the explicit solver.
There are several reasons for this deficiency.
First, the implicit code is substantially newer than our explicit code, and although we have attempted to develop a reasonably well-optimized implementation, there is undoubtedly room for improvement (e.g.~by switching to matrix-free operators where possible).
Moreover, the RAS/RMS smoothers require the use of relatively large subdomains, which results in relatively large computational expenses.
In some cases, the increased robustness of the implicit solvers may still justify their use, as with the implicit code, it is no longer necessary to carefully tune the time step size to avoid instabilities --- a procedure that can substantially increase the time required to set up a complex model.
We anticipate that further work, both in improving the algorithms and their implementations, will make the present scheme useful for low Reynolds number applications.

Treating moderate-to-high Reynolds numbers may require more than simply optimizing our implementation.
One possibility would be to develop an alternative approximation to the Schur complement of Eq.~\eqref{eqn-matrix-impib} for use in the SC smoother.
Alternatively, it may be necessary to reformulate the equations.
For instance, one possibility is that instead of solving Eq.~\eqref{eqn-system-impib}, we instead could solve
\begin{equation}
	\left(\begin{array}{ccc}
		\cA & \cG  & -\cS_h^{n}  \\
		-\cD & \V 0  & \V 0 \\
		-\cJ_h^{n} & \V 0 & \frac{1}{\dt}\cK_h^{-1}
	\end{array}\right)
	\left(\begin{array}{c}
		\u^{n+1}\\
		p^{n+1} \\
		\F^{n+1}
	\end{array}\right) =
	\left(\begin{array}{c}
		\V{g} \\
		0 \\
	 	\V{G}
	\end{array}\right). \label{eqn-system-impib-rev}
\end{equation}
A potential advantage of this formulation is that, for very large stiffnesses, the system is similar to a constrained formulation \cite{BKallemov16-RigidIBAMR, FBalboaUsabiagaXX-rigid_multiblob}.
Effective preconditioners have been developed for this class of problems \cite{BKallemov16-RigidIBAMR, FBalboaUsabiagaXX-rigid_multiblob} and could potentially be extended to the case of FSI with stiff elastic structures.
In the meantime, the development of effective, general-purpose preconditioners for implicit IB formulations with volumetric (codimension-0) structures remains an open problem.

\appendix

\section{Spatial Discretization}

\label{sec-spatial-discretization}

This appendix briefly describes our spatial discretization of the IB equations~\eqref{eqn-momentum}--\eqref{eqn-elastForce}, which is similar to that used in earlier studies \cite{RDGuy15-gmgiib, BEGriffith12-aortic_valve, BEGriffith12-ib_volume_conservation}.

\subsection{Eulerian discretization}
\label{sec-eulerian-discretization}

The Eulerian equations are approximated on a uniform Cartesian grid with grid spacing $h = \dx_1 = \dx_2$ using a staggered-grid discretization in which the Eulerian velocity $\u$ and force $\f$ are approximated at the centers of the Cartesian grid cell edges, 
and the Eulerian pressure is approximated at the centers of the grid cells.
The cell centers are labeled using integer indices $(i,j)$, and the cell edges are labeled using shifted indices, i.e.~$(i-\half,j)$ for $x_1$ edges and $(i,j-\half)$ for $x_2$ edges.
In this notation, $p_{i,j}$ indicates the approximation to $p(\x,t)$ at location $\x_{i,j}$, $(u_1)_{i-\half,j}$ indicates the approximation to the $x_1$ component of the velocity at location $\x_{i-\half,j}$, and $(f_2)_{i,j-\half}$ indicates the approximation to the $x_2$ component of the force at location $\x_{i,j-\half}$.
Spatial Eulerian operators, including the scalar Laplacian $\lap_h$ and vector Laplacian $\V{\grad}^2_h$, gradient $\grad_h$, and divergence $\grad_h \cdot \mbox{}$ are discretized using standard second-order finite differences.
Physical boundary conditions are treated in a manner described previously \cite{BEGriffith09-efficient}.

\subsection{Lagrangian discretization}
\label{sec-lagrangian-discretization}

The Lagrangian force density $\F(\s,t)$ defined in Eq.~\eqref{eqn-F-functional} is discretized on a curvilinear mesh that is free to cut through the background Eulerian grid as the structure moves.
The structure is discretized using a collection of Lagrangian nodes labeled by integer indices $(l,m)$, and we associate to each node curvilinear mesh spacings $(\Delta s_1, \Delta s_2)$.
Simple finite difference approximations are used to evaluate the Lagrangian forces, as described previously~\cite{BEGriffith05-ib_accuracy, BEGriffith07-ibamr_paper, BEGriffith12-aortic_valve, BEGriffith12-ib_volume_conservation}.
Specifically, an approximation to the derivative in the $s_1$ direction of a Lagrangian variable $\Phi(\s,t)$ is defined at a shifted ``half-index'' location by
\begin{equation}
	(D_{s_1} \Phi)_{l+\half,m} = \frac{\Phi_{l+1,m} - \Phi_{l,m}}{\Delta s_1},
\end{equation}
in which $\Phi_{l,m}$ approximates $\Phi(\s,t)$ at curvilinear mesh node $\s_{l,m}$.
Our tests consider only fibers with a zero resting length, for which the fiber tension $T$ and unit tangent vector $\V{\tau}$ are also approximated at shifted locations by
\begin{align}
	T_{l+\half,m} &= \alpha \left\| (D_{s_1} \X)_{l+\half,m} \right\|, \\
	\V{\tau}_{l+\half,m} &= \frac{(D_{s_1} \X)_{l+\half,m}}{\left\| (D_{s_1} \X)_{l+\half,m} \right\|}.
\end{align}
Using these definitions of $D_{s_1}$, $T$, and $\V{\tau}$, we compute an approximation to $F(\s_{l,m}, t)$ via
\begin{equation}
	\F_{l,m} = (D_{s_1}(T \V{\tau}))_{l,m}.
\end{equation}

\subsection{Lagrangian-Eulerian interaction}

Interaction between Lagrangian and Eulerian variables is mediated by integral transforms~\eqref{eqn-F-f} and \eqref{eqn-u-interpolation}.
In the discrete version of the convolution equations, the singular Dirac delta kernel is replaced by a regularized kernel of the form $\delta_h(\x) = \Pi_{i=1}^{d}\delta_h(x_i)$, in which the one-dimensional regularized kernel is $\delta_h(x_i) = \frac{1}{h} \phi(\frac{x_i}{h})$.
In this work, we use Peskin's four-point regularized delta function \cite{Peskin02}, which is defined in terms of the basic kernel function
\begin{equation}
	\phi(r) =  \begin{cases}
		\frac{1}{8}\left( 3 - 2|r| + \sqrt{1 + 4 |r| - 4r^2} \right),	&  0 \leq |r| < 1  \\
        \frac{1}{8}\left( 5 - 2|r| - \sqrt{-7 + 12 |r| - 4r^2} \right),	&  1 \leq |r| < 2 \\
        0,																&  2 \leq |r|.
        \end{cases}
\end{equation}
In two spatial dimensions, a discretized version of the force spreading equation~\eqref{eqn-F-f} is used to obtain the Eulerian force density $\f$ from $\F = (F_1, F_2)$ on the finest level of the locally refined Cartesian grid via
\begin{align}
(f_1)_{i - \half, j} &= \sum_{l,m} (F_1)_{l, m} \, \delta_h ( \x_{i - \half, j} - \X_{l,m}) \, \Delta s_1 \Delta s_2, \\
(f_2)_{i, j - \half} &= \sum_{l,m} (F_2)_{l, m} \, \delta_h ( \x_{i, j - \half} - \X_{l,m}) \, \Delta s_1 \Delta s_2.
\end{align}
Similarly, the Eulerian fluid velocity $\u = (u_1, u_2)$ is interpolated to the curvilinear mesh on the finest grid level to obtain the structural velocity field $\U = (U_1, U_2)$ via
\begin{align}
(U_1)_{l,m} &= \sum_{i, j} (u_1)_{i - \half, j} \, \delta_h ( \x_{i - \half, j} - \X_{l,m}) \Delta x_1 \Delta x_2,  \\
(U_2)_{l,m} &= \sum_{i, j} (u_2)_{i, j - \half} \, \delta_h ( \x_{i, j - \half} - \X_{l,m}) \Delta x_1 \Delta x_2.
\end{align}
As in the continuous equations, we use the shorthand $\f = \cS_h[\X]\,\F$ and $\U = \cJ_h[\X]\,\u$ for these discretized coupling operators.
Moreover, so long as the operators are evaluated using the same structural configuration, $\cS_h[\X] = \cJ_h^*[\X]$ because the same kernel function appears in both of the discretized integrals.

\section*{Acknowledgements}
A.P.S.B.~and R.D.G.~gratefully acknowledge discussions with Gerry Puckett on related solvers used in geodynamics applications.
We also gratefully acknowledge assistance from Barry Smith in profiling and optimizing the performance of the RAS and RMS smoothers used in this work.

\section*{Bibliography}

\bibliography{boyceg-bib/cardiac.bib,boyceg-bib/ib.bib,boyceg-bib/num_analysis.bib}

\end{document}